\documentclass[11pt]{article}

\usepackage{xypic}
\usepackage{amsgen}
\usepackage{amsmath}
\usepackage{amstext}
\usepackage{amsbsy}
\usepackage{amsopn}
\usepackage{amsfonts}
\usepackage{amssymb}
\usepackage{eepic}
\usepackage[pdftex]{graphicx}
\usepackage{epsf}
\usepackage{pstricks}
\xyoption{all}

\def\Box{\square}
\def\edge{\relbar\joinrel\relbar}

\def\tra#1{\smash{\mathop{\mid\kern
-1pt\joinrel\relbar\joinrel\relbar}\limits^{*}_{#1}}}
\def\longtra#1{\smash{\mathop{\mid\kern
-1pt\joinrel\relbar\joinrel\relbar\joinrel\relbar}\limits^{*}_{#1}}}
\def\vlongtra#1{\smash{\mathop{\mid\kern-1pt\joinrel\relbar\joinrel\relbar\joinrel\relbar\joinrel\relbar}\limits^{*}_{#1}}}
\def\vvlongtra#1{\smash{\mathop{\mid\kern
-1pt\joinrel\relbar\joinrel\relbar\joinrel\relbar\joinrel\relbar\joinrel\relbar}\limits^{*}_{#1}}}
\def\vvvlongtra#1{\smash{\mathop{\mid\kern
-1pt\joinrel\relbar\joinrel\relbar\joinrel\relbar\joinrel\relbar\joinrel\relbar\joinrel\relbar}\limits^{*}_{#1}}}
\def\etra#1{\smash{\mathop{\mid\kern
-1pt\joinrel\relbar\joinrel\relbar}\limits_{#1}}}

\def\iff{\Leftrightarrow}
\def\Rw{\Rightarrow}
\def\oo{\overline}

\def\wt{\widetilde}

\def\dwh{\widetilde}

\def\F{{\cal{F}}}

\def\L{{\cal{L}}} 
\def\M{{\cal{M}}}
\newcommand{\N}{{\rm I}\kern-2pt {\rm N}}
\def\SB{\mathbb{SB}}

\def\S{{\cal{S}}}

\def\levi{\mbox{Levi}\,}

\def\nk{\mbox{c-rk}\,}
\def\cmr{\mbox{cm-rk}\,}
\def\crk{\mbox{c-rk}}

\def\diam{\mbox{diam}\,}
\def\lat{\mbox{Lat}\,}

\def\gth{\mbox{gth}\,}
\def\SC{\mbox{SC}}

\def\rk{\mbox{rk}\,}

\def\geo{\mbox{Geo}\,}
\def\matro{\mbox{Mat}\,}
\def\flats{\mbox{Fl}\,}

\def\per{\mbox{Per}\,}

\def\star{\mbox{St}}
\def\cstar{\overline{\mbox{St}}}

\def\max{\mbox{max}}
\def\maxdeg{\mbox{maxdeg}\,}
\def\mindeg{\mbox{mindeg}\,}
\def\min{\mbox{min}}

\def\het{\mbox{ht}\,}

\def\P{{\cal{P}}}

\def\RR{\mathbb{R}}

\def\G{{\cal{G}}}

\def\p{\varphi}

\def\inv{^{-1}}

\def\bi{\begin{itemize}}
\def\ei{\end{itemize}}
\def\beq{\begin{equation}}
\def\eeq{\end{equation}}

\newtheorem{T}{Theorem}[section]
\newcommand{\bt}{\begin{T}}
\newcommand{\et}{\end{T}}
\newcommand{\ftd}{$\square$\end{T}}

\newtheorem{Proposition}[T]{Proposition}
\newcommand{\bp}{\begin{Proposition}}
\newcommand{\ep}{\end{Proposition}}
\newcommand{\fpd}{$\square$\end{Proposition}}

\newtheorem{Definition}[T]{Definition}
\newcommand{\bd}{\begin{Definition}}
\newcommand{\ed}{\end{Definition}}

\newtheorem{Lemma}[T]{Lemma}
\newcommand{\bl}{\begin{Lemma}}
\newcommand{\el}{\end{Lemma}}
\newcommand{\fld}{$\square$\end{Lemma}}

\newtheorem{Corol}[T]{Corollary}
\newcommand{\bc}{\begin{Corol}}
\newcommand{\ec}{\end{Corol}}
\newcommand{\fcd}{$\square$\end{Corol}}

\newtheorem{Result}[T]{Result}
\newcommand{\br}{\begin{Result}}
\newcommand{\er}{\end{Result}}
\newcommand{\frd}{$\square$\end{Result}}

\newtheorem{Remark}[T]{Remark}
\newcommand{\brem}{\begin{Remark}}
\newcommand{\erem}{\end{Remark}}
\newcommand{\fremd}{$\square$\end{Remark}}

\newtheorem{Example}[T]{Example}
\newcommand{\be}{\begin{Example}}
\newcommand{\ee}{\end{Example}}

\newtheorem{Problem}[T]{Problem}
\newcommand{\bq}{\begin{Problem}}
\newcommand{\eq}{\end{Problem}}

\newcommand{\proof}
   {\par\medbreak\noindent{\bf Proof}.\enspace}

\newcommand{\qed}{
$\Box$
\par\bigbreak}

\newlength{\lengtha} 
\newlength{\lengthb} 

\normalbaselines
\def\abstract#1{\par\bigskip
\begingroup\small
\baselineskip=12truept
\begin{center}ABSTRACT\end{center}
\par\medskip\par\noindent
\null\hfill\hbox{\vbox{\hsize=5truein\noindent#1}}
\hfill\null\par\endgroup\par}

\title{A new notion of vertex independence and rank for finite graphs}
\author{{\bf John Rhodes}\\ 
 $ $\\ {\em Department of Mathematics, University of California, Berkeley,}\\ 
{\em California 94720, U.S.A.}\\
{\em email:} rhodes@math.berkeley.edu, BlvdBastille@aol.com\\
$ $\\
{\bf Pedro V. Silva}\\ $ $\\ {\em Centro de
Matem\'{a}tica, Faculdade de Ci\^{e}ncias, Universidade do
Porto,}\\ {\em R. Campo Alegre 687, 4169-007 Porto, Portugal}\\
{\em email:} pvsilva@fc.up.pt} \date{\today}

\begin{document}
\maketitle

\begin{center}\small
\end{center}

\abstract{A new notion of vertex independence and rank for a finite
  graph $G$ is introduced. The independence of vertices is based on
  the boolean independence of columns of a natural boolean matrix
  associated to $G$. Rank is the cardinality of the largest set of
  independent columns. Some basic properties and some more advanced theorems are
  proved. Geometric properties of the graph are related to its rank and
independent sets.} 
 
\tableofcontents

\section{Introduction}

The background and prehistory for this paper goes something like the
following. In 2006 Zur Izakhian \cite{Izh} defined the notion of
independence for columns (rows) of as matrix with coefficients in a
supertropical semiring. Restricting this concept to the superboolean
semiring $\SB$ (see Subsection 2.4), and then to the subset of boolean
matrices (equals matrices with coefficients 0 and 1), we obtain the
notion of independence of columns (rows) of a boolean matrix. This
notion has several equivalent formulations (see Subsection 2.4 of this
paper and references there), one involving permanent, another being
the following: if $M$ is an $m \times n$ boolean matrix, then a subset
$J$ of columns of $M$ is {\em independent} if and only if there exists
a subset $I$ of rows of $M$ with $|I| = |J| = k$ and the $k \times k$
submatrix $M[I,J]$ can be put into upper triangular form (1's on the
diagonal, 0's strictly above it, and 0's or 1's below it) by
independently permuting the rows and columns of $M[I,J]$. 

This is the notion of independence for columns of a boolean matrix we
will use in this paper. 
In 2008 the first author suggested that this
idea would have application in many branches of Mathematics and
especially in Combinatorial Mathematics. In this paper we apply
it to the vertices of a finite graph. For other applications of this
notion to lattices, posets and matroids by Izhakian and the first author, see
\cite{IR1,IR2,IR3}. 

If $M$ is an $m \times n$ boolean matrix with column space $C$,
then the set $\cal{H}$ of independent subsets of $C$ satisfies the
following axioms (see \cite{IR1,IR2}): 
\bi
\item[(H)] $\cal{H}$ is nonempty and closed under taking subsets
  (making it a {\em hereditary collection});
\item[(PR)] for all nonempty $J, \{ p \} \in {\cal{H}}$, there exists
  some $x \in J$ such that $( J \setminus \{x\}) \cup  \{ p \} \in
  {\cal{H}}$ (the {\em point replacement} property).
\ei
Hereditary collections arising from some boolean matrix $M$ as above
are said to be {\em boolean representable}. 
A very interesting
question is which hereditary collections have boolean representations,
a question which the authors will address in a near future paper
\cite{RSil}. The elementary properties of such boolean representable
collections were considered in \cite{IR1,IR2,IR3} and it was shown in
\cite{IR2}  
that all matroids have boolean representations.

In this paper we restrict our atention to finite graphs (with no loops
and no multiple edges), see Subsection 2.2. However, there are several
ways to define such a graph by a boolean matrix. The one chosen in
matroid theory by Whitney \cite{Whi} and related to the Levi graph is to attach
the boolean matrix $M(G)$ to the graph $G = (V,E)$, where $V$ is the
set of vertices and $E$ the set of edges considered as 2-sets of $V$,
with $M(G)$ the $|V| \times |E|$ boolean matrix defined by $M(G)(v,e)
= 1$ if $v$ lies in $e$, and 0 otherwise. Now whether we consider the
columns of $M(G)$ as independent in our boolean sense or in the
usual vector space sense (over the field 
$\mathbb{Z}_2$), we obtain the same independent sets which form a
matroid called a graphical matroid, see \cite{IR1,Oxl,Oxl2}. 

So this viewpoint has been extensively worked out \cite{Oxl, Oxl2},
and mainly following Tutte's suggestions, ideas from graphs like
connectedness ($n$-connected) can be extended to matroids, etc.

A perhaps more obvious way to associate a boolean matrix to a graph $G
= (V,E)$ is via the $|V| \times |V|$ boolean adjacency matrix (see
Subsection 2.2) $A_G = (a_{ij})$, where $a_{ij} = 1$ if $\{ i,j\}$ is
an edge of $G$, and 0 otherwise. So $A_G$ can be an arbitrary 
symmetric square boolean matrix with 0's in the main
diagonal (see also \cite{BT}). However, in this paper we choose
$A_G^c$ which is $A_G$ with 
0 and 1 interchanged. This approach is indicated from the lattice/poset
case \cite{IR3}, and that finite boolean modules (equals semilattices) have dual
spaces which separate points and the dual space is reversing the
order, see \cite[Chapter 9.1 and 9.2]{RS}. 

Also if $A_G$ were used, then $K_n$ (the complete graph on $n$
vertices) and its complement $\oo{K_n}$ would have sets of 2 or less
vertices being the independent sets or only the empty set being
independent respectively, clearly not a good choice.

Thus our new notion of independence of a subset of vertices $X
\subseteq V$ of a graph $G = (V,E)$ is that the columns corresponding
to $X$ in $A_G^c$ are boolean independent. Note that, by using
$A_G^c$, all subsets of vertices of $K_n$ are independent. This is
termed {\em c-independence} for vertices of $G$, and the cardinality
of the largest independent set of vertices ia termed {\em c-rank},
denoted $\nk$. Note that we work with the superboolean
semiring $\SB$, for representation by matrices over $GF(2)$ the reader
can be referred to a recent paper by Brijder and Traldi \cite{BT}. 

As we mentioned before, Whitney associated to each finite graph $G =
(V,E)$ a (graphical) matroid \cite{Whi}. In this paper we more or less
reverse this procedure and treat each graph $G$ as given ``like a
matroid'' in the following manner. The graph $G$ has the boolean
representation $A_G^c = M$. Each boolean representation $M$, see
\cite{IR2, IR3, RSil}, gives rise
to the {\em lattice of flats} (see Subsection 2.2) of $M$. This
corresponds to the idea in matroid theory of the geometric lattice of
flats of a matroid (see \cite{Oxl}). Given the boolean matrix $M$ with
column space $C$, the lattice of flats of $M$ consists of the subsets
of $C$ corresponding to where the rows of $M$ are zero, closed under
all intersections (see Subsection 2.2). Then Theorem \ref{indhei} yields
that the independent subsets of $C$ with respect to $M$ are the
partial transversals of the partition of successive differences for
some maximal chain of the lattice of flats. This relates to earlier
ideas of  Bjorner and Ziegler \cite{BZ}.  

If $L$ is the geometric lattice of
flats of a matroid $P = (C,{\cal{H}})$, then taking the boolean
representation $M_L$ corresponding to $L$ and restricted to the atom
generators ($M_L$ is $A_L^c$ -- where $A_L$ is the $L$ incidence
matrix -- restricted to the atom rows $C$, then transposed so
considered as columns), then the lattice of flats of the boolean
matrix $M_L$ as described before is the same as the geometric lattice of
usual flats of the matroid (see \cite{IR2}). Thus this approach truly 
generalizes the matroid approach.  

When applying this ``boolean combinatorics'' approach to some standard
field of Mathematics (e.g. finite graph theory), usually the notion of
rank is well known, the notion of independence is new, and the
approach tends quickly to some well developed subfield of the subject
under study. Somehow geometry is also supposed to show up in this approach:
see below!

Enough of the general background. In this paper the boolean
representation for a graph is $A_G^c$ and the notions of
c-independence and c-rank are taken with respect to $A_G^c$. The
lattice of flats for the graph $G = (V,E)$ can be realised by closing
$\{ \star(v) \mid v \in V \}$ under all intersections, where $\star(v)
= \{ v' \in V \mid \{ v, v'\} \in E \}$. This and other preliminaries
are done in Section 2. The c-rank and how to calculate the
c-independent subsets of vertices are discussed in Section 3. It is
proved that the c-rank is the height of the lattice of flats and
c-independent subsets can be calculated by Theorem
\ref{indhei}(iv)-(v).

In Section 4 we characterize graphs of low c-rank.
 Section 5 is devoted to the interesting case of sober
connected graphs of c-rank 3 (we call a graph {\em sober} if the
mapping $\star$ is injective). 
In Section 6, our new notions acquire a distinctive geometric flavor
in connection with Levi graphs and partial euclidean geometries. 
$\geo$ is defined in the appropriate context and $\geo$ of the
Petersen graph is computed to be the Desargues configuration, see
Example \ref{petdes}. Section 7 collects results concerning
cubic graphs, including characterizations of the graphs whose lattice
of flats satisfies the most famous lattice-theoretic
properties.  A variation of the
concept of c-rank appropriate to deal with minors is discussed in
Section 8. Finally, Section 9 relates a graph and its complement
graph in the context of our new notions.

\section{Preliminaries}

\subsection{Posets and lattices}

Our lattice and poset terminology is more or less standard (see
\cite{GHK,Gra,MMT,RS}).  For ease of exposition we assume all
posets, lattices and graphs to be finite, although many of the results 
admit generalizations to the infinite case.

Given a finite poset $(P,\leq)$ and $p,q \in P$, we say that $p$ {\em
  covers} $q$ if $p > q$ but there is no $r \in P$ such that $p > r >
q$. It is standard to represent finite posets by means of their {\em
  Hasse diagram}: in this directed graph, the vertices are the
elements of $P$ and $(p,q)$ is an edge when $p$ covers $q$. Note that
a chain in $(P,\leq)$ is maximal if and only if it corresponds to some
path in the Hasse diagram connecting a maximal element to a minimal element.  

The {\em height} of $(P,\leq)$ is defined by
$$\het P = \max\{ k \in \N \mid p_0 < p_1 < \ldots < p_k \mbox{ is a
  chain in }P \}.$$
Equivalently, $\het P$ is the maximum length of a path in the Hasse
diagram of $P$.  

We say that  $(P,\leq)$ is a lattice if, for all $p,q \in P$, there
exist
$$\begin{array}{l}
p\vee q = \min\{ x \in P \mid x \geq p,q \},\\
p\wedge q = \max\{ x \in P \mid x \leq p,q \}.
\end{array}$$
If only the first (respectively the second) of these conditions is
satisfied, we talk of a $\vee$-{\em semilattice} (respectively
$\wedge$-{\em semilattice}). 
We say that $P' \subseteq P$ constitutes a {\em sublattice} of
$(P,\leq)$ if $p\vee q, p\wedge q \in P'$ for all $p,q \in P'$. 
Note there need be no relation between the top (bottom) of $P'$ and of
$P$.  Every point is a sublattice.  

A lattice $(L,\leq)$ is said to be {\em distributive} if
$$p \wedge (q \vee r) = (p \wedge r) \vee (p \wedge r)$$
holds for all $p,q,r \in L$, a condition
which is equivalent to its dual. We shall
say that  $(L,\leq)$ is {\em 
  modular} if there is no sublattice of the form
$$
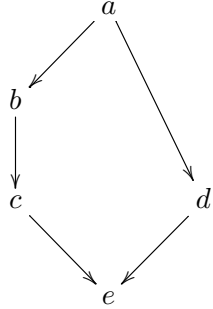
\xymatrix{
& a \ar[dl] \ar[ddr] & \\
b \ar[d] && \\
c \ar[dr] && d \ar[dl] \\
& e &
}$$
If we only exclude such sublattices when $d$ covers $e$, the lattice is {\em
  semimodular}. 
It is well known that every distributive lattice is modular, and
modular implies of course semimodular.

An {\em atom} of $L$ is an element covering the minimum element
$0$. A semimodular latice is called {\em geometric} if every element
is a join of atoms (0 being the join of the empty set). Finally, $L$
satisfies the {\em Jordan-Dedekind condition} if all the maximal
chains in $L$ have the same length.

\subsection{Graphs}
\label{graphs}

Throughout this paper, graphs are finite, undirected, and have neither
loops nor multiple edges. Formally, a (finite) graph is an ordered pair $G =
(V,E)$, where $V$ is a (finite) set (the set of vertices) and $E
\subseteq \{ X \in 2^V: |X| = 2\}$ (the set of edges). In other words,
the edges are 2-{\em subsets} of $V$ (an $n$-subset is a subset with
$n$ elements). We assume the
reader to be familiar with the basic concepts of graph theory (see
e.g. \cite{Die}).

Clearly, $(2^V, \subseteq)$ is a distributive lattice with $X \wedge Y = X
\cap Y$ and $X \vee Y = X
\cup Y$. If $X \subseteq V$ has $k$ elements, we say it is a $k$-{\em
  subset} of $V$.

Given $\S \subseteq 2^V$,
it is easy to see that
$$\widehat{\S} = \{ \cap S \mid S \subseteq \S \}$$
is the $\wedge$-subsemilattice of $(2^V, \subseteq)$ generated by
$\S$. Note that $\cap \S = \min \widehat{\S}$, and also $V =
\cap \emptyset = \max \widehat{\S}$. In fact, $(\widehat{\S},\subseteq)$ is
itself a lattice with
$$P \vee Q = \cap \{ X \in \S \mid P\cup Q \subseteq X\}.$$
However,  $(\widehat{\S},\subseteq)$ is not in general a sublattice of
$(2^V, \subseteq)$ since $P \vee Q$ (in $(\widehat{\S},\subseteq)$) needs
not be $P \cup Q$ (see \cite{GHK,RS}). 

Note that
\beq
\label{bht}
\het\widehat{\S} \leq |\S|
\eeq
since any chain in $\widehat{\S}$ is necessarily of the form
$$V \supseteq X_1 \supset X_1\cap X_2 \supset \ldots \supset X_1 \cap
\ldots \cap X_k$$
for distinct $X_1,\ldots,X_k \in \S$.

Finally, we say that $\{ y_1,\ldots,y_k \}$ is a {\em 
  transversal}  of the partition of the successive differences for the
chain $X_0 \supset \ldots \supset X_k$ in 
$\widehat{\S}$ if $y_i \in X_{i-1} \setminus X_i$ for $i =
1,\ldots,k$. A subset of a transversal is a {\em 
partial transversal}.

Given $v \in V$, the {\em star} of $v$ is defined by
$$\star(v) = \{ w \in V \mid w \mbox{ is adjacent to $v$ in } G\}.$$
More generally, given $W \subseteq V$, we write 
$$\star(W) = \displaystyle\cap_{w\in W} \star(w).$$
Note that $\star(\emptyset) = V$. Let $\S_W = \{ \star(w) \mid w \in W
\}$. It is immediate that
\beq
\label{hatstar}
\widehat{\S_W} = \{ \star(W') \mid W' \subseteq W\}.
\eeq 
We call a subset of the form $\star(W)$ $(W \in V)$ a {\em flat} and
say that $\widehat{\S_V}$ is the {\em lattice of flats} of $G$, also
denoted by $\flats G$.  We
believe this to be a new concept for graphs.  

Note that, for a connected graph $G = (V,E)$, we can define a metric
$d$ on $V$ by
$$d(v,w) = \mbox{ length of the shortest path connecting $v$ and $w$
  (counting edges)}.$$
The {\em diameter} of $G$, denoted by $\diam G$, is the maximum value
in the image of $d$.

Given a
finite graph $G$, the {\em girth} of $G$, denoted by $\gth G$, is the
length of the shortest cycle in $G$ (assumed to be $\infty$ is $G$ is
acyclic). Note that $\gth G \geq 3$ for any finite graph.

We shall use the notation 
$$\hat{n} = \{ 1,\ldots,n\}$$
throughout the paper.
Assume now that $V = \hat{n}$. The {\em adjacency matrix} of
$G = (V,E)$ is the $n\times n$ boolean matrix $A_G =(a_{ij})$ defined
by 
$$a_{ij} = \left\{
\begin{array}{ll}
1&\mbox{ if }\{i,j\} \in E\\
0&\mbox{ otherwise}
\end{array}
\right.$$
The matrix $A_G^c$ is obtained by interchanging 0 and 1 all over
$A_G$. If the graph is clear from the context, we shall write just $A$
and $A^c$.

\subsection{Matroids}
\label{matroi}

Let $V$ be a set and let $X \subseteq 2^V$. We say that $X$ is a {\em
  hereditary collection} if $X$ is closed under taking subsets. The
hereditary collection is said to be a {\em matroid} if the following
condition (the {\em 
  exchange property}) holds:
\bi
\item[(EP)] For all $I,J \in X$ with $|I| = |J|+1$, there exists some
  $i \in I\setminus J$ such that $J \cup \{ i \} \in X$.
\ei

There are many other equivalent definitions of matroid. For details,
the reader is referred to \cite{Oxl}.

\subsection{Superboolean matrices}

Following \cite{IR1}, we shall view boolean matrices as matrices over
the {\em superboolean semiring} $\SB = \{ 0,1,1^{\nu} \}$, where
addition and multiplication are described respectively by 
$$\begin{tabular}{l|lll}
+&0&1&$1^{\nu}$\\
\hline
0&0&1&$1^{\nu}$\\
1&1&$1^{\nu}$&$1^{\nu}$\\
$1^{\nu}$&$1^{\nu}$&$1^{\nu}$&$1^{\nu}$
\end{tabular}
\hspace{1.5cm}
\begin{tabular}{l|lll}
$\cdot$&0&1&$1^{\nu}$\\
\hline
0&0&0&0\\
1&0&1&$1^{\nu}$\\
$1^{\nu}$&0&$1^{\nu}$&$1^{\nu}$
\end{tabular}$$
We denote by $\M_n(\SB)$ the set of all $n \times n$ matrices with
entries in $\SB$.  Note that $n \times n$ boolean matrices are {\em not}
a subsemiring of  $\M_n(\SB)$ since $1 + 1  = 1^{\nu}$.

Next we present definitions of independency and rank appropriate to
the context of superboolean matrices, introduced in \cite{Izh} (see
also \cite{IR1}).  

We say that vectors $C_1,\ldots,C_m \in \SB^n$ are {\em dependent} if
$\lambda_1C_1+ \ldots \lambda_mC_m \in \{ 0,1^{\nu} \}$ for some
$\lambda_1,\ldots, \lambda_m \in \{0,1\}$ not all zero. Otherwise,
they are said to be {\em independent}.   

Let $S_n$ denote the symmetric group on $\hat{n}$. The {\em
  permanent} of a matrix $M \in \M_n(\SB)$ (a positive version of the
determinant) is defined by 
$$\per M = \displaystyle\sum_{\sigma \in S_n} \prod_{i=1}^n
m_{i,i\sigma}.$$
Recall that addition and
multiplication take place in the semiring $\SB$ defined above.

Given $I,J \subseteq \hat{n}$, we denote by $M[I,J]$ the submatrix of
$M$ with entries $m_{ij}$ $(i\in I, j \in J)$. In particular,
$M[\hat{n},j]$ denotes the $j$th column vector of $M$ for each $j \in \hat{n}$.

\bp
\label{nons}
{\rm \cite[Th. 2.10]{Izh}, \cite[Lemma 3.2]{IR1}}
The following conditions are equivalent for every $M \in \M_n(\SB)$:
\bi
\item[(i)] the column vectors $M[\hat{n},1],\ldots, M[\hat{n},n]$ are
  independent; 
\item[(ii)] ${\rm Per}\,M = 1$;
\item[(iii)] $M$ can be transformed into some lower triangular matrix
  of the form 
\beq
\label{nons1}
\left(
\begin{matrix}
1&&0&&0&&\ldots&&0\\
?&&1&&0&&\ldots&&0\\
?&&?&&1&&\ldots&&0\\
\vdots&&\vdots&&\vdots&&\ddots&&\vdots\\
?&&?&&?&&\ldots&&1
\end{matrix}
\right)
\eeq
by permuting rows and permuting columns independently.
\ei
\ep

A square matrix satisfying the above (equivalent) conditions is said to be
{\em nonsingular}. 

Given (equipotent) $I,J \subseteq \hat{n}$, we say that $I$ is a {\em
  witness} for $J$ in $M$ if $M[I,J]$ is nonsingular. 

\bp
\label{indwit}
{\rm \cite[Th. 3.11]{Izh}}
The following conditions are equivalent for all $M \in \M_n(\SB)$ and
$J \subseteq \{ 1,\ldots,n\}$: 
\bi
\item[(i)] the column vectors $M[\hat{n},j]$ $(j \in J)$ are independent;
\item[(ii)] $J$ has a witness in $M$.
\ei
\ep

The subsets of independent column vectors of a given superboolean
matrix, which include the empty subset and are closed for subsets,
constitute an important example of a hereditary
  collection. Hereditary  collections which have boolean
  representations, which include matroids as a 
very important particular case, were discussed in \cite{IR1,IR2,IR3}
and will be also the object of a future paper by the present authors,
seeking necessary and sufficient conditions.

\bp
\label{altrank}
{\rm \cite[Th. 3.11]{Izh}}
The following are equal for a given $M \in \M_n(\SB)$:
\bi
\item[(i)] the maximum number of independent column vectors in $M$;
\item[(ii)] the maximum number of independent row vectors in $M$;
\item[(iii)] the maximum size of a subset $J \subseteq \hat{n}$ having
  a witness in $M$;  
\item[(iv)] the maximum size of a nonsingular submatrix of $M$.
\ei
\ep

The {\em rank} of a matrix $M \in \M_n(\SB)$, denoted by $\rk M$, is
the number described above. A row of $M$ is called an $n$-{\em marker}
if it has one entry 1 and all the remaining entries are 0. The
following remark follows from Proposition \ref{nons}: 

\bc
\label{mark}
{\rm \cite[Cor. 3.4]{IR1}}
If $M \in \M_n(\SB)$ is nonsingular, then it has an $n$-marker.
\ec

\section{The c-rank of a graph}

In this section, we assume that $G = (V,E)$ denotes a finite graph
with $V = \hat{n}$. 

The following result prepares the ground for an important connection
between matrix rank and the height of the lattice of flats as defined
in Subsection \ref{graphs}, and will acquire great relevance in the
study of independence. This relates to earlier work by Bjorner and
Ziegler \cite{BZ}.

\bt
\label{indhei}
Given a finite graph $G$, the
following conditions are equivalent for every $J \subseteq \hat{n}$: 
\bi
\item[(i)] the column vectors $A^c[\hat{n},j]$ $(j \in J)$ are independent;
\item[(ii)] $J$ has a witness in $A^c$;
\item[(iii)] ${\rm ht}\, \widehat{\S_J} = |J|$;
\item[(iv)] $J$ is a transversal of the partition of successive
  differences for some chain of {\rm Fl}$\, G$; 
\item[(v)] $J$ is a partial transversal of the partition of successive
  differences for some maximal chain of {\rm Fl}$\, G$.
\ei
\et

\proof
We may assume that $J$ is nonempty.

(i) $\iff$ (ii). By Proposition \ref{indwit}.

(ii) $\Rw$ (iii). Let $I$ be a witness for $J$ in $A^c$. Permuting
rows and columns if necessary, we may assume that 
$A^c[I,J]$ is of the form (\ref{nons1}). Write $k = |J|$ and let
$i_1,\ldots,i_k$ and $j_1,\ldots,j_k$ indicate the new ordering of
rows and columns in the reordered matrix. Then the reordered $A[I,J]$ is of
the form 
\beq
\label{nons2}
\left(
\begin{matrix}
0&&1&&1&&\ldots&&1\\
?&&0&&1&&\ldots&&1\\
?&&?&&0&&\ldots&&1\\
\vdots&&\vdots&&\vdots&&\ddots&&\vdots\\
?&&?&&?&&\ldots&&0
\end{matrix}
\right)
\eeq
and so
$$\star(j_r,\ldots,j_k) \cap I = \{ i_1, \ldots, i_{r-1} \}$$
for $r = 1,\ldots,k$. 
Since $j_k \notin \star{j_k}$, we get
\beq
\label{indhei1}
V = \star(\emptyset) \supset \star{j_k} \supset \star(j_{k-1}, j_k)
\supset \ldots \supset \star(j_1,\ldots,j_k) 
\eeq
and so  $\het \widehat{\S_J} \geq k$. Hence  $\het \widehat{\S_J} = k$
by (\ref{bht}). 

(iii) $\Rw$ (ii). In view of (\ref{bht}), it is easy to see that we
must have necessarily a chain of the form (\ref{indhei1}), where $J =
\{ j_1,\ldots, j_k\}$ and the $j_r$ are all distinct. For $r =
1,\ldots,k-1$, take $i_r \in \star(j_{r+1},\ldots,j_k) \setminus
\star(j_r)$, and also $i_k = j_k$. With the rows (respectively
columns) ordered by $i_1,\ldots,i_k$ (respectively $j_1,\ldots,j_k$),
the matrix $A^c[I,J]$ is now of the form (\ref{nons1}) and so $I$ is a
witness for $J$ in $A^c$. 

(ii) $\Rw$ (iv). Let $I = \{ i_1, \ldots, i_k \}$ be a witness for
$J = \{ i_1, \ldots, j_k \}$ in $A^c$. Similarly to
the proof of (ii) $\Rw$ (iii), we may assume that $A[I,J]$ is of
the form (\ref{nons2})
and so
$$\star(i_1,\ldots,i_r) \cap J = \{ j_{r+1}, \ldots, j_{k} \}$$
for $r = 1,\ldots,k$. Hence
\beq
\label{indhei3}
V = \star(\emptyset) \supset \star(i_1) \supset \star(i_1, i_2)
\supset \ldots \supset \star(i_1,\ldots,i_k) 
\eeq
is a chain in $\flats G$. 
Since $j_r \in
\star(i_1,\ldots,i_{r-1})\setminus \star(i_r)$, then $J$ is a 
transversal for (\ref{indhei3}).

(iv) $\Rw$ (ii). Assume that $J = \{ i_1, \ldots, j_k \}$ is a
transversal for a chain 
$$\star(X_0) \supset \star(X_1) \supset \ldots \supset \star(X_k)$$
 in $\flats G$. We may assume that $j_r \in \star(X_{r-1}) \setminus
 \star(X_r)$ for $r = 1,\ldots,k$. Then, for each $r$, there exists
 $i_r \in X_r$ such that $j_r \notin \star(i_r)$. However, if $s < r$,
 then $j_r \in \star(X_{r-1}) \subseteq \star(X_s) \subseteq
 \star(i_s)$ and it follows easily that, with the rows (respectively
columns) ordered by $i_1,\ldots,i_k$ (respectively $j_1,\ldots,j_k$),
the matrix $A^c[I,J]$ is now of the form (\ref{nons1}) and so $I$ is a
witness for $J$ in $A^c$.

(iv) $\Rw$ (v). Since a partial transversal of a maximal chain is a
transversal for some subchain of the original chain.

(v) $\Rw$ (iv). Since every chain can be refined to get a maximal
chain.
\qed

To simplify terminology, we say that the vertices $j_1, \ldots, j_k
\in \hat{n}$ are {\em c-independent} if the column vectors
$A^c[\hat{n},j_1], \ldots, A^c[\hat{n},j_k]$ are independent. 

\brem
\label{oneortwo}
Let $G$ be a finite graph and let $j_1,j_2 \in \hat{n}$. Then:
\bi
\item[(i)] $j_1$ is c-independent;
\item[(ii)] $j_1,j_2$ are c-independent if and only if {\rm St}$(j_1)
  \neq$ {\rm St}$(j_2)$. 
\ei
In particular,  $j_1,j_2$ are c-independent if they are adjacent.
\erem

\proof
(i) This follows from every column vector in $A^c$ being nonzero due
to the absence of loops in $G$. 

(ii) Since every column vector in $A^c$ is nonzero, it follows from
Theorem \ref{indhei} that $j_1,j_2$ are c-independent if and only if
$A^c[\hat{n},j_1]$ and $A^c[\hat{n},j_2]$ are distinct,
i.e. $\star(j_1) \neq \star(j_2)$. 
\qed

\bt
\label{rkht}
Let $G = (V,E)$ be a finite graph. Then {\rm rk}$\, A^c =$ {\rm
  ht}$\,${\rm Fl}$\, G$.  
\et

\proof
Let $k = \rk A^c$. Then there
exists some $J \subseteq \hat{n}$ such that $|J| = k$ and the column vectors
$A^c[\hat{n},j]$ $(j \in J)$ are independent. Hence 
$\het \widehat{\S_J} = k$ by Theorem \ref{indhei}. Since $\het
\widehat{\S_J} \leq \het  \widehat{\S_V}$ by (\ref{hatstar}), it
follows that  $\rk A^c \leq \het \flats G$. 

Assume now that $\het \flats G = \ell$. Then there exists a
(maximal) chain 
\beq
\label{rkht1}
V = \star(\emptyset) \supset \star(V_{\ell}) \supset \ldots \supset \star(V_1)
\eeq
for some $V_1,\ldots, V_{\ell} \subseteq V$. We claim that there exist
$j_1,\ldots, j_{\ell} \in V$ such that 
\beq
\label{rkht2}
\star(j_r, \ldots, j_{\ell}) =
\star(V_r)
\eeq
for $r = 1,\ldots, \ell$. Indeed, since $\star(V_{r+1}) \supset
\star(V_r)$, we can take $j_r \in V_r$ such that
$\star(V_{r+1}) \not\subseteq \star(j_r)$. 
Writing $V_{\ell +1}
= \emptyset$, we proceed by induction on $r = \ell, \ldots,1$:
assume that (\ref{rkht2}) holds for $r+1$. Hence 
$$\star(V_{r+1}) = \star(j_{r+1}, \ldots, j_{\ell}) \supset \star(j_r,
\ldots, j_{\ell}) = \star(j_r) \cap \star(j_{r+1}, \ldots, j_{\ell})
\supseteq \star(V_r)$$ 
and so $\star(V_r) = \star(j_r, \ldots, j_{\ell})$ by the maximality
of the length of the chain (\ref{rkht1}). Thus (\ref{rkht2}) holds.

Take $J = \{ j_1,\ldots,j_{\ell} \}$. Since  $\het
\widehat{\S_J} \leq \het
\widehat{\S_V} = \ell$, it follows from (\ref{rkht1}) and
(\ref{rkht2}) that $\het 
\widehat{\S_J} = \ell = |J|$ and so 
 the column vectors $A^c[\hat{n},j]$ $(j \in J)$ are independent by
 Theorem \ref{indhei}. Thus $\het \flats G = \ell \leq \rk A^c$
 and so $\rk A^c = \het \flats G$.  
\qed

We say that the above number is the {\em c-rank} of the graph $G$ and
we denote it by $\nk G$. Note that, in view of Theorem \ref{indhei},
$\nk G$ is also the maximum size of a (partial) transversal of the
partition of successive differences of a (maximal) chain of $\flats G$.

We present now some straightforward properties of the c-rank of a graph.
Let $\maxdeg G$ (respectively $\mindeg G$) denote the maximum
(respectively minimum) degree of a vertex in $G$.

\bp
\label{maxd}
Let $G$ be a finite graph. Then {\rm c}-{\rm rk}$\, G \leq {\rm
  maxdeg}\, G +1$. 
\ep

\proof
Since in a chain of the form (\ref{rkht1}), we have necessarily
$|\star(V_{\ell})| \leq \maxdeg G$. 
\qed

\bp
\label{cc}
Let $G$ be a finite graph with connected components
$G_1,\ldots,G_m$. Then {\rm c}-{\rm rk}$\, G = \max\{$ {\rm c}-{\rm
  rk}$\, G_1, \; \ldots,$ {\rm c}-{\rm rk}$\,G_m \; \}$. 
\ep

\proof
Since in any chain of the form (\ref{rkht1}), the $V_r$ and the
$\star(V_r)$ must necessarily be taken in one same connected
component. 
\qed

In view of this result, we may focus our attention, from now on, on
{\em connected} graphs.
 
We say that $G' = (V',E')$ is a {\em subgraph} of $G = (V,E)$ if $V'
\subseteq V$ and $E' \subseteq E$ (up to isomorphism!). If  $V'
\subseteq V$ and  $E' = E 
\cap 2^{V'}$, we say that $G'$ is a {\em restriction} of $G$. 

Given graphs $G = (V,E)$ and  $G' = (V',E')$, a {\em morphism} $\p:G
\to G'$ is a mapping $\p:V \to V'$ such that $v\p \edge w\p$ is an
edge of $G'$ whenever $v \edge w$ is an
edge of $G$. We say that $\p$ is a {\em retraction} if $G'$ is a
restriction of $G$ and $\p|_{V'}$ is the identity mapping.

\bp
\label{res}
Let $G = (V,E)$, $G' = (V',E')$ be finite graphs. 
\bi
\item[(i)] If $G'$ is a restriction of $G$, then {\rm c}-{\rm rk}$\,G' \leq$
  {\rm c}-{\rm rk}$\, G$. 
\item[(ii)] If $G'$ is a complete subgraph of $G$, then {\rm c}-{\rm rk}$\,
  G \geq$ {\rm c}-{\rm rk}$\, G' = |V'|$. 
\ei
\ep

\proof
(i) If $G'$ is a restriction of $G$, then any (nonsingular) submatrix
of $A_{G'}^c$ is also a (nonsingular) submatrix of $A_{G}^c$. 

(ii) A complete subgraph of $G$ is necessarily a restriction, hence
$\nk G \geq \nk G'$ by part (i). The equality $\nk G' = |V'|$ follows
from the following fact: if $K_n$ denotes the complete graph on $n$
vertices, then $A_{K_n}^c$ is the identity matrix. 
\qed

Note that  $\nk G' \leq \nk G$ may not hold if $G'$
is a mere subgraph of $G$. For instance, it is easy to check that the
square
$$
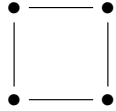
\xymatrix{
\bullet \ar@{-}[r] \ar@{-}[d] & \bullet \ar@{-}[d] \\
\bullet \ar@{-}[r] & \bullet
}$$
has c-rank 2, but after removing an
edge the c-rank increases (cf. Proposition \ref{casetwo}).

We introduce now a concept that will ease the discussion of c-rank in
many circumstances. We call a finite graph $G = (V,E)$ {\em sober} if
the star mapping $\star: V \to 2^v$ is injective. The following remark
is immediate from Remark \ref{oneortwo}:

\brem
\label{consim}
The following conditions are equivalent for a finite connected graph $G$: 
\bi
\item[(i)] $G$ is sober;
\item[(ii)] all 2-subsets of vertices of $G$ are independent.
\ei
\erem

\bp
\label{soberrest}
Let $G = (V,E)$ be a finite connected graph. Then $G$ admits a
retraction onto a sober
connected restriction $G' = (V',E')$ such that {\rm Fl}$\,G \cong$
{\rm Fl}$\,G'$. 
\ep

\proof
Let $V'$ be a cross-section for the star mapping $\star:V \to 2^V$
of $G$ and
let $G'$ be the restriction of $G$ induced by $V'$. It is
straightforward that $G'$  
is isomorphic to the graph having as vertices the
equivalence classes of $V$ induced by $\star$ 
and edges $X \edge Y$ whenever $x \edge y$ is an edge of $G$ for some
$x \in X$ and $y \in Y$.

For every $v \in
V$, let $v' \in V'$ be the unique vertex in $V'$ such that $\star(v')
= \star(v)$. We claim that, for all $v,w \in V$,
\beq
\label{edges}
\{ v,w \} \in E \iff \{ v',w' \} \in E'.
\eeq
Indeed, if $v 
\edge w$ is an edge in $G$, then so is $v
\edge w'$ and therefore $v' \edge w'$. 

Conversely, assume that $\{ v',w' \} \in E' \subseteq E$. Then we
successively get $\{ v,w' \} \in E$ and $\{ v,w \} \in E$, hence
(\ref{edges}) holds. Thus $\p: V \to V'$ defined by $v\p = v'$ is a
graph morphism from $G$ to the restriction $G'$, indeed a retraction.

Moreover, any path
$v_1 \edge \ldots \edge v_n$ in $G$ induces a path $v'_1 \edge
\ldots \edge v'_n$ in $G'$ and so $G'$ is connected.

Let $\star':V' \to 2^{V'}$ denote the star mapping
of $G$. Suppose that $v,w \in V$ are such that $\star'(v') =
\star'(w')$. It follows from (\ref{edges}) that
$$\star(v) = \{ z \in V \mid z' \in \star'(v') \} = \{ z \in V \mid z'
\in \star'(w') \} = \star(w),$$
hence $v' = w'$ and so $G'$ is sober.

We claim that
$$\begin{array}{rcl}
\theta:(\flats G,\subseteq)&\to&(\flats G',\subseteq)\\
\star(W)&\mapsto&\star'(W')
\end{array}$$
is an isomorphism of posets (and therefore of lattices).

It is immediate that $\theta$ is surjective and preserves order. It
remains to show that $\theta$ is well defined and injective.

For every $W \subseteq V$, it follows from (\ref{edges}) that
$$\star(W) = \displaystyle\cap_{w\in W} \star(w) = \cap_{w\in W}\{ x
\in V \mid x' \in \star'(w') \} = \{ x
\in V \mid x' \in \star'(W') \},$$
$$\star'(W') =  \{ x \in V \mid x' \in \displaystyle\cap_{w\in W} \star'(w') \}'
 =  \{ x \in V \mid x \in \displaystyle\cap_{w\in W} \star(w) \}' =
 (\star(W))'.$$ Therefore $\star(W) = \star(Z) \iff \star'(W') =
 \star'(Z')$ holds for all $W,Z \subseteq V$ and so $\theta$ is an
 isomorphism.
\qed 

However, the restriction in Proposition \ref{soberrest} does not need
to be unique (up to 
isomorphism). For instance, the graph
$$
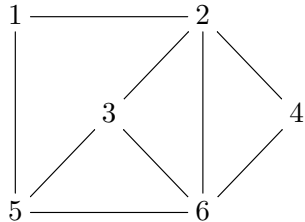
\xymatrix{
1 \ar@{-}[rr] \ar@{-}[dd] && 2 \ar@{-}[dl] \ar@{-}[dr] \ar@{-}[dd] & \\
& 3 \ar@{-}[dl] \ar@{-}[dr] && 4 \ar@{-}[dl] \\
5 \ar@{-}[rr] && 6 &
}$$
is itself sober and connected (and has $\mindeg$ 2), and so it is the
restriction obtained  
by removing vertex 1. It is easy to check that the star lattices of
both graphs are isomorphic and of the form:
$$
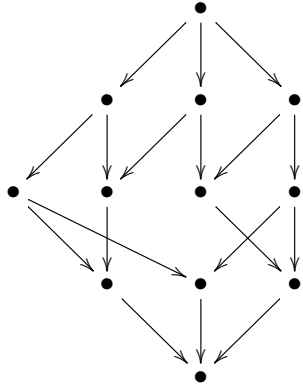
\xymatrix{
&& \bullet \ar[dl] \ar[d] \ar[dr] & \\
& \bullet \ar[dl] \ar[d] & \bullet  \ar[dl] \ar[d] & \bullet  \ar[dl]
\ar[d] \\
\bullet \ar[dr]  \ar[drr] & \bullet \ar[d] & \bullet \ar[dr] & \bullet \ar[dl]
\ar[d] \\
& \bullet \ar[dr] & \bullet \ar[d] & \bullet \ar[dl] \\
&& \bullet &
}$$

It is easy to characterize sober trees. Recall that a vertex of
degree 1 is called a {\em leaf}.

\bp
\label{strees}
A tree $T = (V,E)$ is sober if and only if no two leafs are at
distance 2 from each other.
\ep

\proof
Indeed, assume that $v,w \in V$ are distinct. If $\star(v) = \star(w)$
and has more than one element, then $T$ would admit a square and would
not be a tree, hence $\star(v) = \star(w)$ can only occur if both $v$
and $w$ are leafs, in which case $\star(v) = \star(w)$ is equivalent
to $d(v,w) = 2$.
\qed

We establish next an inductive relation that may prove useful in the
computation of the c-rank. Given a graph $G = (V,E)$, and
 $J \subseteq V$, write $\crk_G J = \rk
A^c[V,J]$. Note that $\nk G = \crk_G V$. 

We recall also that, for $X \subseteq V$, the graph $G-X$ is obtained
from $G$ by removing all the vertices in $X$ and all the edges adjacent
to them.

\bt
\label{inwit}
Let $G = (V,E)$ be a finite graph and $m \geq
2$. Then the following conditions are equivalent: 
\bi
\item[(i)] {\rm c}-{\rm rk}$\, G \geq m$.
\item[(ii)] There exist $v,w \in V$ such that:
\bi
\item
{\rm St}$(v) \neq$ {\rm St}$(w)$;
\item
{\rm c}-{\rm rk}$_{G-\{v,w\}}(${\rm St}$(v) \cap$ {\rm St}$(w)) \geq m-2$.
\ei
\ei
\et

\proof
(i) $\Rw$ (ii). If $\nk G \geq m$, then by Proposition \ref{indwit} there
exist $I,J \subseteq V$ 
such that $A^c[I,J]$ is nonsingular and $|J| =
m$. In view of Proposition \ref{nons}, we may reorder the rows
(respectively columns) of $A^c[I,J]$ by $i_1,\ldots,i_m$
(respectively $j_1,\ldots,j_m$) to get a matrix of the form
(\ref{nons1}). Since $j_2 \in \star(i_1) \setminus \star(i_2)$, we
have $\star(i_1) \neq \star(i_2)$. On the other hand, $i_1,i_2 \notin
\{ i_3, \ldots,i_m\} \cup \{ j_3, \ldots, j_m\}$ since $i_1,i_2 \in
\star( j_3, \ldots, j_m)$, hence $\{  i_3, \ldots,i_m\}$ is a witness
for $\{ j_3, \ldots, j_m\}$ in $G\setminus \{ i_1,i_2\}$. Therefore  
$j_3,\ldots,j_m \in \star(i_1)
\cap \star(i_2)$ are c-independent in $G\setminus \{ i_1,i_2\}$ and so
condition (ii) holds.  

(ii) $\Rw$ (i). Since $\crk_{G-\{v,w\}}(\star(v) \cap \star(w)) \geq m-2$, there
exist distinct $j_3,\ldots,j_m \in \star(v) \cap \star(w)$ and
$i_3,\ldots,i_m \in V\setminus\{ v,w \}$ such that
$A^c[i_3,\ldots,i_m;j_3,\ldots,j_m]$ 
is nonsingular. Reordering rows and columns if necessary, we may
assume that $A^c[i_3,\ldots,i_m;j_3,\ldots,j_m]$ is of the form
(\ref{nons1}). Since $\star(v) \neq \star(w)$, we may assume that there exists
some $j_2 \in \star(v) 
\setminus \star(w)$ and take $i_1 = j_1 = v$ and $i_2 = w$. It is
straightforward to check that $A^c[i_1,\ldots,i_m;j_1,\ldots,j_m]$ is
of the form (\ref{nons1}), hence condition (i) holds for $J = \{ j_1,
\ldots, j_m \}$. 
\qed

\section{Low c-rank}

We start analyzing the sober cases and go as far as characterizing
c-rank 4. In view of Propositions \ref{cc} and \ref{soberrest}, in the
discussion of c-rank $\geq 3$ we pay special attention to the case of
sober connected graphs.  
 
\bp
\label{trivials} 
Let $G = (V,E)$ be a finite graph. Then:
\bi
\item[(i)] {\rm c}-{\rm rk}$\, G = 0$ if and only if $V = \emptyset$.
\item[(ii)] {\rm c}-{\rm rk}$\, G = 1$ if and only if $V \neq \emptyset$ and $E
  = \emptyset$. 
\ei
\ep

\proof
Clearly, $\nk G \geq 0$ under all circumstances and the empty graph
has c-rank 0. On the other hand, if $V \neq \emptyset$, then $A^c$ has
at least one 1 in the diagonal, yielding $\nk G \geq 1$. This proves
(i). 
Moreover, if $E \neq \emptyset$, it follows from Remark
\ref{oneortwo}(ii) that $\nk G \geq 2$, thus (ii) holds as well. 
\qed

We recall that a graph $G = (V,E)$ is called {\em bipartite} if $V$
admits a nontrivial partition $V = V_1 \cup V_2$ such that  
$$E \subseteq \{ \{ v_1,v_2\} \mid v_1 \in V_1,v_2 \in V_2 \}.$$
If this inclusion can be made to be an equality, the graph is said to
be {\em complete bipartite}.  

\bp
\label{casetwo} 
Let $G = (V,E)$ be a finite graph with $E \neq \emptyset$. Then the
following conditions are equivalent: 
\bi
\item[(i)] {\rm c}-{\rm rk}$\, G = 2$;
\item[(ii)] $G$ has no subgraph $v_1 \edge v_2 \edge v_3$ with
  {\rm St}$(v_1) \neq$ {\rm St}$(v_3)$; 
\item[(iii)] $G$ is a disjoint union of complete bipartite graphs;
\item[(iv)] $G$ has no restrictions of the following forms:
$$
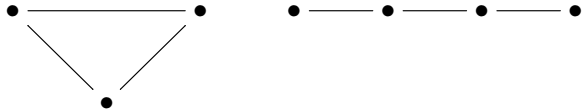
\xymatrix{
\bullet \ar@{-}[rr] \ar@{-}[dr] && \bullet & \bullet \ar@{-}[r] &
\bullet \ar@{-}[r] & \bullet \ar@{-}[r] & \bullet \\ 
& \bullet \ar@{-}[ur] &&&&&  
}$$
\ei
\ep

\proof
(i) $\Rw$ (ii). Suppose that $G$ has a subgraph $v_1 \edge v_2 \edge v_3$ with
$\star(v_1) \neq \star(v_3)$. We may assume that there exists some $w \in \star(v_1) \setminus \star(v_3)$. Consider the chain
$$V = \star(\emptyset) \supset \star(v_1) \supset \star(v_1,v_3) \supset \star(v_1,v_2,v_3)$$
(taking respectively $v_1, w, v_2$ to show that the inclusions are strict). Thus $\nk G \geq 3$.

(ii) $\Rw$ (iii). Since (ii) holds, any path
$$v_1 \edge v_2 \edge \ldots \edge v_k$$
must satisfy $\star(v_1) = \star(v_3) = \ldots$ and $\star(v_2) = \star(v_4) = \ldots$
Since we may also assume $G$ to be connected, it follows that we can take a pair of adjacent edges $(w_1,w_2)$ and partition $V = V_1 \cup V_2$ by
$$V_1 = \{ v \in V \mid \star(v) = \star(w_1)\},\quad
V_2 = \{ v \in V \mid \star(v) = \star(w_2)\}.$$

Since $w_2 \in \star(w_1)$, we have $w_2 \in \star(v)$ for every $v \in V_1$. Hence $V_1 \subseteq \star(w_2)$ and so $V_1 \subseteq \star(v)$ for every $v \in V_2$. On the other hand, if $v,v' \in V_1$ are adjacent, then $v \in \star(v') = \star(v)$, a contradiction. Similarly, no two vertices in $V_2$ can be adjacent. Thus $G$ is complete bipartite.

(iii) $\Rw$ (iv). It is well-known that no bipartite graph admits cycles of odd length. Suppose that $G = (V,E)$ is bipartite complete (with respect to the partition $V = V_1 \cup V_2$) and has a restriction of the form 
\beq
\label{casetwo1}
v_1 \edge v_2 \edge v_3 \edge v_4.
\eeq
Then $\{ v_1,v_3\}$ and  $\{ v_2,v_4\}$ belong to the different sides of the partition and so there exists an edge $v_1 \edge v_4$ in $G$, contradicting (\ref{casetwo1}) being a restriction. Therefore $G$ can have no restriction of the form (\ref{casetwo1}) either.

(iv) $\Rw$ (i). Suppose that $\nk G \geq 3$. After reordering, $A^c$ has a submatrix of the form
$$\begin{tabular}{l|lll}
$i_1$&1&0&0\\
$i_2$&?&1&0\\
$i_3$&?&?&1\\
\hline
&$j_1$&$j_2$&$j_3$
\end{tabular}$$
and so  $A$ has a submatrix of the form
$$\begin{tabular}{l|lll}
$i_1$&0&1&1\\
$i_2$&?&0&1\\
$i_3$&?&?&0\\
\hline
&$j_1$&$j_2$&$j_3$
\end{tabular}$$
Then $j_2 \edge i_1 \edge j_3$ is a subgraph with 3 distinct
vertices. We may assume that the triangle $K_3$ is not a restriction
of $G$. Since $i_2 \edge j_3$ is an edge, it follows that $i_2 \neq
j_2$. Hence  
\beq
\label{casetwo2}
j_2 \edge i_1 \edge j_3 \edge i_2
\eeq
is a subgraph of $G$ with 4 distinct vertices. Since there is no edge
$i_2 \edge j_2$ and $K_3$ is not a restriction of $G$, then
(\ref{casetwo2}) is a restriction of $G$ and so (iv) fails as required.
\qed

\bp
\label{gfour} 
Let $G$ be a finite graph. Then the following conditions are equivalent:
\bi
\item[(i)] {\rm c}-{\rm rk}$\, G \geq 4$;
\item[(ii)] $G$ has a subgraph 
$$
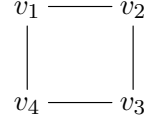
\xymatrix{
v_1 \ar@{-}[r] \ar@{-}[d] & v_2 \ar@{-}[d] \\
v_4 \ar@{-}[r] & v_3
}$$
with {\rm St}$(v_1) \neq$ {\rm St}$(v_3)$ and {\rm St}$(v_2) \neq$
{\rm St}$(v_4)$. 
\ei
\ep

\proof
Write $G = (V,E)$.

(i) $\Rw$ (ii). Suppose that $\nk G \geq 4$. After reordering, $A$ has a submatrix of the form
$$\begin{tabular}{l|llll}
$i_1$&0&1&1&1\\
$i_2$&?&0&1&1\\
$i_3$&?&?&0&1\\
$i_4$&?&?&?&0\\
\hline
&$j_1$&$j_2$&$j_3$&$j_4$
\end{tabular}$$
Then we have edges
$$
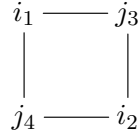
\xymatrix{
i_1 \ar@{-}[r] \ar@{-}[d] & j_3 \ar@{-}[d] \\
j_4 \ar@{-}[r] & i_2
}$$
in $G$. Since $i_3 \in
\star(j_4) \setminus \star(j_3)$ and $j_2 \in \star(i_1) \setminus
\star(i_2)$, the vertices $i_1,i_2,j_3,j_4$ are all distinct and (ii) holds. 

(ii) $\Rw$ (i). If (ii) holds, then we may assume out of symmetry that
there exist some $w \in \star(v_1) \setminus \star(v_3)$ and $z \in
\star(v_2) \setminus \star(v_4)$.   Consider the chain
$$V \supset \star(v_1) \supset \star(v_1,v_3)
\supset \star(v_1,v_3,z) \supset \emptyset$$
(taking respectively $v_1, w, v_4, v_2$ to show that the
inclusions are strict). Thus $\nk G \geq 4$. 
\qed

Now Propositions \ref{casetwo} and \ref{gfour} combined provide a
characterization of c-rank 3.

\bp
\label{gfive} 
Let $G$ be a finite graph. Then the following conditions are equivalent:
\bi
\item[(i)] {\rm c}-{\rm rk}$\, G \geq 5$;
\item[(ii)] $G$ has a subgraph 
$$
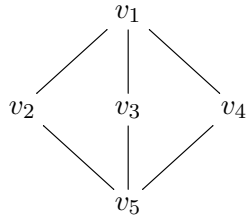
\xymatrix{
&v_1 \ar@{-}[dl] \ar@{-}[d] \ar@{-}[dr] & \\
v_2 \ar@{-}[dr] & v_3 \ar@{-}[d] & v_4 \ar@{-}[dl] \\
& v_5 &
}$$
with {\rm St}$(v_1) \neq$ {\rm St}$(v_5)$, {\rm St}$(v_2) \neq$ {\rm
  St}$(v_3)$ and {\rm St}$(v_2) \;\cap$ {\rm St}$(v_3) \not\subseteq$
{\rm St}$(v_4)$. 
\ei
\ep

\proof
Write $G = (V,E)$.

(i) $\Rw$ (ii). 
Suppose that $\nk G \geq 5$. After reordering, $A$ has a submatrix of the form
$$\begin{tabular}{l|lllll}
$i_1$&0&1&1&1&1\\
$i_2$&?&0&1&1&1\\
$i_3$&?&?&0&1&1\\
$i_4$&?&?&?&0&1\\
$i_5$&?&?&?&?&0\\
\hline
&$j_1$&$j_2$&$j_3$&$j_4$&$j_5$
\end{tabular}$$
Then 
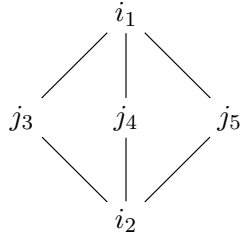
$$\xymatrix{
&i_1 \ar@{-}[dl] \ar@{-}[d] \ar@{-}[dr] & \\
j_3 \ar@{-}[dr] & j_4 \ar@{-}[d] & j_5 \ar@{-}[dl] \\
& i_2 &
}$$
is a subgraph of $G$ with 5 distinct
vertices. Now $j_2 \in
\star(i_1) \setminus \star(i_2)$, $i_4 \in \star(j_5) \setminus
\star(j_4)$ and $i_3 \in (\star(j_4) \cap \star(j_5)) \setminus 
\star(j_3)$ and so (ii) holds. 

(ii) $\Rw$ (i). If (ii) holds, then we may assume out of symmetry that
there exist some $w_1 \in \star(v_1) \setminus \star(v_5)$, $w_2 \in
\star(v_2) \setminus \star(v_3)$ and $w_3 \in
(\star(v_2) \cap \star(v_3)) \setminus \star(v_4)$. Consider the chain
$$V \supset \star(v_1) \supset \star(v_1,v_5)
\supset \star(v_1,v_5,w_3) \supset \star(v_1,v_5,w_3,w_2) \supset \emptyset$$
(taking respectively $v_1, w_1, v_4, v_3, v_2$ to show that the
inclusions are strict). Thus $\nk G \geq 5$. 
\qed

Now Propositions \ref{gfour} and \ref{gfive} combined provide a
characterization of c-rank 4.

We can now use the previous results to give a complete
characterization of sober connected graphs with low c-rank (in view of
(iv), see \cite{CFS}):

\bc
\label{sclow} 
Let $G = (V,E)$ be a finite sober connected graph. Then:
\bi
\item[(i)] {\rm c}-{\rm rk}$\, G = 0$ if and only if $G$ is the empty graph;
\item[(ii)]  {\rm c}-{\rm rk}$\, G = 1$ if and only if $G \cong K_1$;
\item[(iii)]  {\rm c}-{\rm rk}$\, G = 2$ if and only if $G \cong K_2$;
\item[(iv)]  {\rm c}-{\rm rk}$\, G = 3$ if and only if $|E| \geq 2$
  and $G$ has no squares;
\item[(v)]  {\rm c}-{\rm rk}$\, G = 4$ if and only if $G$ has a square
  but no subgraph 
$$
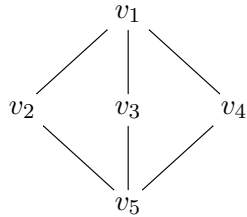
\xymatrix{
&v_1 \ar@{-}[dl] \ar@{-}[d] \ar@{-}[dr] & \\
v_2 \ar@{-}[dr] & v_3 \ar@{-}[d] & v_4 \ar@{-}[dl] \\
& v_5 &
}$$
with {\rm St}$(v_2)\; \cap$ {\rm St}$(v_3) \not\subseteq$ {\rm St}$(v_4)$.
\item[(vi)]  {\rm c}-{\rm rk}$\, G \geq 5$ if and only if $G$ has a
  subgraph of the above form. 
\ei
\ec

\proof
(i) and (ii) follow immediately from Proposition \ref{trivials}. Since
Since sober connected nontrivial complete bipartite graphs
can have only one edge, (iii) follows from Proposition \ref{casetwo}.

Now part (iii) implies that $\nk G \geq 3$ if and only if $|E| \geq
2$, and so (iv) follows from Proposition \ref{gfour}.

Finally, Propositions \ref{gfour} and \ref{gfive} yield (v) and (vi).
\qed

\section{The c-independent subsets in c-rank 3}
\label{gvx}

We shall denote by $\SC n$ the class of all finite sober connected
graphs of c-rank $n$.
Throughout this section, all graphs are in $\SC 3$. In view of
Corollary \ref{sclow}(iv), these graphs have no squares (for such
graphs with few vertices, see
\cite{CFS}). 

The following
lemma collects some elementary facts involving this class of
graphs. We recall that a graph is called {\em cubic} if all vertices
have degree 3. 

\bl
\label{sct}
Let $G$ be a finite connected graph.
\bi
\item[(i)]
If $G$ is cubic and  {\rm gth}$\, G \geq 5$, then $G \in$  {\rm SC}$\, 3$.
\item[(ii)]
If $G = (V,E) \in$  {\rm SC}$\, 3$, then $|${\rm St}$(v)\; \cap$
 {\rm St}$(w)| \leq 1$ holds for all distinct vertices $v,w$ of $G$.
\ei
\el

\proof
(i) If $G$ is non sober, then $G$ would contain a square,
contradicting $\gth G \geq 5$. Hence $G$ is sober. The claim now
follows from Corollary \ref{sclow}.

(ii) Suppose that $|\star(v,w)| > 1$ for distinct vertices $v,w$ of $G
= (V,E)$. Since $G$ 
is sober, we may assume that $\star(v,w) \subset
\star(v)$. Let $a,b \in \star(v,w)$ be distinct. Since $G$ 
is sober, we may assume that there exists some $c \in \star(b)
\setminus \star(a)$. Hence
$$\star(v) \supset \star(v,w) \supset \star(v,w,c) \supset
\star(v,w,c,a)$$
is a chain in $\flats G$, contradicting $\nk G = 3$. 
\qed

By c-rank, the c-independent subsets of a graph $G = 
(V,E)$ in $\SC 3$ can have at
most 3 elements. However, as it will become clear soon enough, the
c-independent subsets of $V$ do not constitute a matroid. Our first
result associates a matroid to $G$: we define $\matro G$ to contain:
\bi
\item
all the $i$-subsets of $V$ for $i \leq 2$;
\item
all the 3-subsets $W$ of $V$ such that
$$\forall v \in V\; W \not\subseteq \star(v).$$
\ei
Note that the latter condition is
equivalent to $\star(W) = \emptyset$. 

\bp
\label{matro}
Let $G \in$  {\rm SC}$\, 3$. Then {\rm Mat}$\, G$ is a matroid.
\ep

\proof
Let $I,J \in \matro G$. Without loss of generality, we may assume that $|I|
= 3$ and $|J| = 2$. Write $I = \{ i_1, i_2, i_3\}$ and $J = \{
j_1,j_2\}$. Suppose that $\{
j_1,j_2,i_k\} \notin \matro G$ for $k = 1,2,3$. Then there exists some $v_k
\in \star(j_1,j_2,i_k) \subseteq \star(j_1,j_2)$. By Lemma \ref{sct}(ii), we get
$v_1 = v_2 = v_3$ and so $I \subseteq \star(v_1)$, contradicting 
$I \in \matro G$. Therefore $\{ j_1,j_2,i_k\} \in \matro G$ for some $k
\in \hat{3}$ and so $\matro G$ is a matroid.
\qed

We identify next the c-independent subsets of
vertices for graphs in $\SC 3$. We say that a 3-subset $P
\subseteq V$ is a {\em potential 
  line} if $|P \cap \star(v)| \leq 1$ for every $v \in V$.

\bt
\label{idind}
Let $G = (V,E)$ be a graph in {\rm SC}$\, 3$ and let $W
\subseteq V$. Then the following conditions are equivalent:
\bi
\item[(i)] $W$ is c-independent;
\item[(ii)] $|W| \leq 2$ or
$$|W| = 3,\; {\rm St}(W) = \emptyset \mbox{ and $W$ is not a potential line}.$$ 
\ei
\et

\proof 
Since $G$ is sober, and by Remark \ref{oneortwo}, $W$ is
c-independent if $|W| \leq 2$. On the other hand, since $\nk G = 3$,
then $V$ has no c-independent 4-subsets. Therefore we may assume that
$|W| = 3$. Write $W = \{ w_1,w_2,w_3\}$. 

Assume that $W$ is independent. By Theorem \ref{indhei}, we may assume
that 
\beq
\label{idind1}
V \supset \star(w_1) \supset \star(w_1,w_2) \supset
\star(w_1,w_2,w_3)
\eeq
is a chain in $\flats G$. Since  $\star(w_1,w_2,w_3) \neq \emptyset$
would allow us to adjoin the empty set to the chain and contradict
c-rank 3, then $\star(W) = \emptyset$. On the other hand, for $v \in 
\star(w_1,w_2)$, we get $|W \cap \star(v)| \geq 2$ and so $W$ is not a
potential line either.

Conversely, assume that $\star(W) = \emptyset$ and $W$ is not a potential
line. Then $|W \cap \star(v)| \geq 2$ for some in $v \in V$. We may
assume that $w_1,w_2 \in \star(v)$. Furthermore, since $G$ is sober,
we may also assume that $\star(w_1) \supset \star(w_1,w_2)$.
(\ref{idind1}) is a chain in $\flats G$ and so $W$ is independent by
Theorem \ref{indhei}. 
\qed

Now Proposition \ref{matro} and Theorem \ref{idind} yield:

\bc
\label{idindc}
Let $G = (V,E)$ be a graph in {\rm SC}$\, 3$. If $G$ has no potential
lines, then the set of all c-independent 
subsets of $V$ constitutes a matroid.
\ec

In view of this result, it is only natural to enquire which graphs
in the above class have no potential lines. It turns out that
diameter makes the difference: 

\bp
\label{nopot}
Let $G$ be a graph in {\rm SC}$\, 3$.
\bi
\item[(i)] If {\rm diam}$\, G < 3$, then $G$ has no potential lines.
\item[(ii)] If {\rm diam}$\, G > 5$, then $G$ has potential lines.
\item[(iii)] If {\rm diam}$\, G \in \{ 3,4,5\}$, then both cases may occur.
\ei
\ep

\proof
(i) First, we note that if $P \subseteq V$ is a potential line and $p,q
\in P$ are distinct, then $d(p,q) \neq 2$ (if $p \edge v \edge q$ is a
path in $G$, then $|P \cap \star(v)| \geq 2$), and if $d(p,q) = 1$,
then the edge $p \edge q$ can lie in no triangle. Hence, if $\diam G <
3$ and $P = \{ p_1,p_2,p_3\}$ is a potential line, then $d(p_1,p_2) =
d(p_1,p_3) = d(p_2,p_3) = 1$ immediately gets us into a
contradiction. Thus (i) holds.

(ii) Assume that $\diam G > 5$. Then $G$ has a {\em geodesic} (path of
minimum length connecting the extreme vertices) of length 6, say
$$v_0 \edge v_1 \edge v_2 \edge v_3 \edge v_4 \edge v_5 \edge v_6$$
Since $d(v_0,v_3) = d(v_3,v_6) = 3$, it follows that $\{
v_0,v_3,v_6\}$ is a potential line of $G$.

(iii) It is enough to show that there exist in $\SC 3$
a graph $G_3$ with diameter 3
and potential lines, and a graph $G_5$ with diameter 5
and no potential lines.

We can take $G_3$ to be the cubic graph
$$
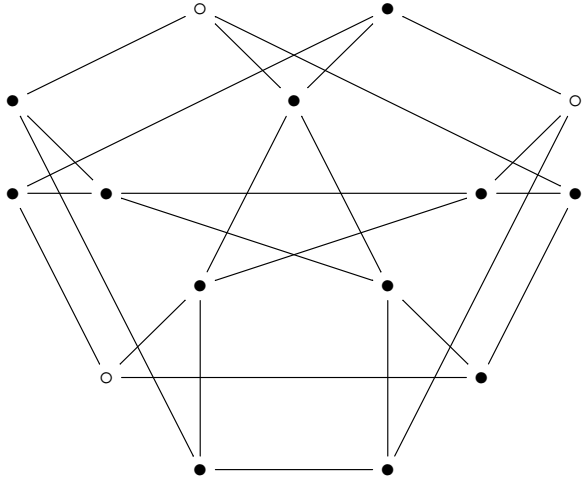
\xymatrix{
&& \circ \ar@{-}[dll] \ar@{-}[dr] \ar@{-}[ddrrrr] && \bullet \ar@{-}[ddllll]
\ar@{-}[dl] \ar@{-}[drr] && \\
\bullet \ar@{-}[ddddrr] \ar@{-}[dr] &&& \bullet \ar@{-}[ddl] \ar@{-}[ddr] &&&
\circ \ar@{-}[dl] \ar@{-}[ddddll] \\
\bullet \ar@{-}[r] \ar@{-}[ddr] & \bullet \ar@{-}[rrrr] \ar@{-}[drrr] &&&&
\bullet \ar@{-}[dlll] \ar@{-}[r] & \bullet \ar@{-}[ddl] \\
&& \bullet \ar@{-}[dl] \ar@{-}[dd] && \bullet \ar@{-}[dd] \ar@{-}[dr] && \\
& \circ \ar@{-}[rrrr] &&&& \bullet & \\
&& \bullet \ar@{-}[rr] && \bullet &&
}$$
Since $\gth G = 5$, it follows from Lemma \ref{sct}(ii) that $G \in
\SC 3$. Straightforward checking shows that $\diam G = 3$ and $G$ has
potential lines such  
as the one defined by the hollow circles.

On the other hand, we can take $G_5$ to be the graph
$$
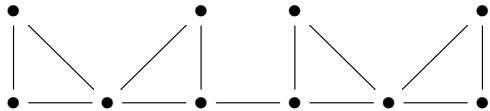
\xymatrix{
\bullet \ar@{-}[d] \ar@{-}[dr] && \bullet \ar@{-}[d] \ar@{-}[dl] &
\bullet \ar@{-}[d] \ar@{-}[dr] && \bullet \ar@{-}[d] \ar@{-}[dl] \\
\bullet \ar@{-}[r] & \bullet \ar@{-}[r] & \bullet \ar@{-}[r] & \bullet
\ar@{-}[r] & \bullet \ar@{-}[r] & \bullet
}$$
It follows easily from Corollary \ref{sclow}(iv) that $G \in
\SC 3$, and it is immediate that $\diam G = 5$. Suppose that $G$ has a
potential line $P$. Then at least two points of $P$ would have to fit
into a subgraph of the form  
$$
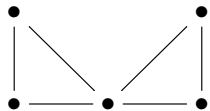
\xymatrix{
\bullet \ar@{-}[d] \ar@{-}[dr] && \bullet \ar@{-}[d] \ar@{-}[dl] \\
\bullet \ar@{-}[r] & \bullet \ar@{-}[r] & \bullet
}$$
leading at once to a contradiction. Therefore $G$ has no potential
lines as claimed.
\qed

\be
\label{coimbra}
Let $G = (V,E)$ be the graph
$$\xymatrix{
1 \ar@{-}[dd] \ar@{-}[rr] \ar@{-}[dr] && 4 \ar@{-}[dr] & \\
& 3 \ar@{-}[dl] \ar@{-}[r] & 5 \ar@{-}[r] & 7 \ar@{-}[dl] \\
2 \ar@{-}[rr] && 6 &
}$$
\ee

\noindent
(see \cite{CFS}). By Corollary \ref{sclow}(iv), we have $G \in \SC
3$. The lattice of 
flats of $G$ can be depicted as 
$$
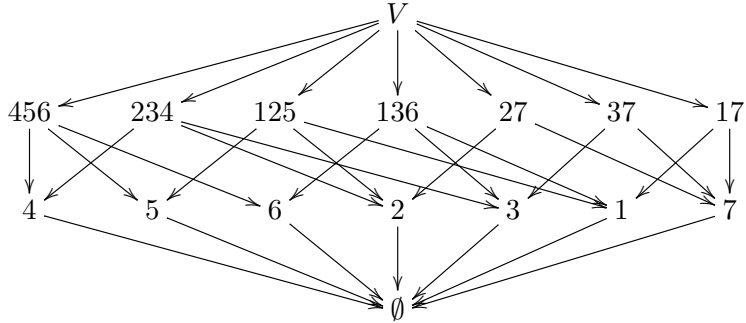
\xymatrix{
&&& V \ar[dlll] \ar[dll] \ar[dl] \ar[d] \ar[dr] \ar[drr] \ar[drrr] &&& \\ 
456 \ar[d] \ar[dr] \ar[drr] & 234 \ar[drr] \ar[drrr] \ar[dl] & 125 \ar[dl]
\ar[dr] \ar[drrr] & 136 \ar[dl] \ar[dr] \ar[drr] & 27 \ar[dl] \ar[drr] & 37
\ar[dl] \ar[dr] & 17 \ar[dl] \ar[d] \\
4 \ar[drrr] & 5 \ar[drr] & 6 \ar[dr] & 2 \ar[d] & 3 \ar[dl] & 1
\ar[dll] & 7 \ar[dlll] \\ 
&&& \emptyset &&&
}$$
It is straightforward to check that $G$ has no potential lines and
$\matro G$ contains all the $i$-subsets of $V$ for $i \leq 3$ except
the flats 125, 136, 234 and 456. In view of Theorem \ref{idind}, these
are precisely the c-independent subsets of $V$. See further remarks
after Corollary \ref{leviconf} relating to the Fano plane. 

\medskip

If we restrict our attention to cubic graphs, the range is a bit
reduced. A list of all cubic graphs up to 12 vertices can be found in
\cite{W6}, where the handy LCF notation is explained and used.

\bc
\label{nopotcub}
Let $G$ be a cubic graph in {\rm SC}$\,3$.
\bi
\item[(i)] If {\rm diam}$\, < 3$, then $G$ has no potential lines.
\item[(ii)] If {\rm diam}$\, G > 3$, then $G$ has potential lines.
\item[(iii)] If {\rm diam}$\, G = 3$, then both cases may occur.
\ei
\ec

\proof
(i) By Proposition \ref{nopot}(i).

(ii) Suppose now that $\diam G > 3$. Let $a,b \in V$ be such that $d(a,b) =
4$, and write $\star(a) = \{ x,y,z \}$. Clearly, $b$ is at distance
$\geq 3$ from $x$, $y$ or $z$. To prevent $\{ a,b,x\}$ from being a
potential line, $a \edge x$ must lie in some triangle. If we try to
avoid other potential lines, also $a \edge y$ and $a \edge z$ must lie
in triangles. Now it is easy to see that at least two of the vertices
$x,y,z$ must be connected through edges. Without loss of generality,
we may assume that $$
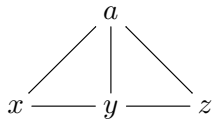
\xymatrix{
& a \ar@{-}[dr] \ar@{-}[d] \ar@{-}[dl] & \\
x \ar@{-}[r] & y \ar@{-}[r] & z
}$$
is a subgraph of $G$. But then we have a square in a sober graph,
contradicting c-rank 3 in view of Corollary \ref{sclow}(iv). Thus (ii)
holds.  

(iii) The example $G_3$ in the proof of Proposition \ref{nopot}(iii)
is cubic, belongs to $\SC 3$, has diameter 3 and has potential lines. 

However, the Heawood graph \cite{W1} is cubic, bipartite, has diameter 3 and
girth 6 (and so is in $\SC 3$, see Proposition \ref{cubic}(iii) in
next section). Suppose that $P = \{ a,b,c\}$ is a potential line of the
Heawood graph. Then the distance between any two distinct vertices in
$P$ cannot be 2, and so must be 1 or 3 in view of the diameter being
3. Thus we obtain a cycle of odd length in the graph, contradicting
the fact of being bipartite. Therefore the Heawood graph has no
potential lines.
\qed

\section{The Levi graph and partial euclidean geometries}

Given a finite graph $G = (V,E)$ we can consider $V$ as ``points'' and
$E$ as ``lines'', where $v$ is on $e$ $(v \in e)$ if and only
if $e$ is incident to $v$, and so $(V,E)$ gives {\em some sort of
  geometry} (see \cite{Cam,Cox}). So the Levi viewpoint for ``lines'' in a
graph is different from our view of taking $\star(v)$ as lines. I this
section, we benefit from this other approach and introduce right away
the concept of partial euclidean geometry.

Let $P$ be a finite nonempty
set and let $\L$ be a nonempty subset of $2^P$. We shall always assume
that $P \cap 2^P = \emptyset$. We say that $(P,\L)$
is a {\em partial euclidean geometry} (abbreviated to PEG) if
the following axioms are satisfied:
\bi
\item[(G1)] $P \subseteq \cup \L$;
\item[(G2)] if $L,L' \in \L$ are distinct, then $|L \cap L'| \leq 1$;
\item[(G3)] $|L| \geq 2$ for every $L \in \L$.
\ei 
The elements of $P$ are called {\em points} and the elements of $\L$
are called {\em lines}.  Given $p \in P$, we denote by $\L(p)$ the set
of all lines 
containing $p$.

The concept of PEG is  an abstract
combinatorial generalization of the following geometric situation: 

Consider a finite set of lines $\L$ in the euclidean space
$\RR^n$. Consider also a finite subset $P$ of $\cup\L 
\subset \RR^n$ such that:
\bi
\item if $L,L' \in \L$ are distinct, then $|L \cap L'| \leq 1$;
\item if $L,L' \in \L$ and $L \cap L' = \{ p \}$, then $p \in P$;
\item $|L \cap P| \geq 2$ for every $L \in \L$.
\ei
Representing each $L\in \L$ by $L\cap P$, it follows that $(\L,P)$ constitutes a
PEG. It is well known that not all PEG's can be represented over an
euclidean space (nor 
any field) (see \cite[Section 2.6]{Gru}).

Using Coxeter's notation (see \cite{Cox}), we say that the PEG $(P,\L)$ is an
$(m_c,n_d)$ {\em configuration} if:
\bi
\item
there are $m$ points and $n$ lines;
\item
each point belongs to $c$ lines;
\item
each line contains $d$ points.
\ei
Hence $cm = dn$, which equals the number of 1's in the (boolean) {\em incidence
  matrix} of $(P,\L)$, where rows are labelled by points and columns by
lines. 

An important example is provided by the famous {\em Desargues
  configuration}. A simple way of defining it is by taking points as
2-subsets of $\hat{5}$ and lines as 3-subsets of $\hat{5}$
(identifying $\{ a,b,c\}$ with $\{ \{ a,b \}, \{ a,c \}, \{ b,c \}
\}$). For a geometric representation, see e.g. \cite{W4}. It is clear
that the Desargues configuration is a $(10_3,10_3)$ configuration. It
has many interesting properties, such as being self-dual (by
exchanging points and lines, we get an isomorphic configuration), and
the automorphism group acts transitively on both vertices and
edges. And it is of course related to the famous Desargues' Theorem
\cite{W4}. Notice that, for every point $p$, there are exactly 3
points noncolinear with $p$ (i.e., not belonging to some line
simultaneously with $p$), and that these 3 points constitute a line!

Now, for every $G = (V,E) \in \SC 3$, let $$\L_G = \{ W \in \flats G
\setminus \{ V \} : |W|
\geq 2 \}$$ and let $\geo G = (V,\L_G)$.

\bp
\label{peg}
If $G \in$ {\rm SC}$3$, then {\rm Geo}$\, G$ is a PEG.
\ep

\proof
Let $v \in V$. Since $\nk G = 3$, then $\star(v) \neq \emptyset$. If
all the elements of $\star(v)$ have degree 1, then 
$G$ sober implies that $v$ has also degree 1 and so $G \cong K_2$,
contradicting Corollary \ref{sclow}(iv). Hence there exists some $w
\in \star(v)$ with degree $\geq 2$ and so $v \in \star(w) \in
L_G$. Thus $\geo G$ satisfies axiom (G1) (and also $\L_G \neq
\emptyset$).

Finally, (G2) follows from Lemma \ref{sct}(ii) and (G3) holds
trivially. Therefore $\geo G$ is a PEG.
\qed

Note that $\matro G$ consists of all subsets of $V$ with at most 2
elements plus all 3-subsets which are contained in no line of $\geo G$.

\bc
\label{cuco}
If $G \in$ {\rm SC}$3$ is cubic with $n$ vertices, then {\rm Geo}$\, G$ is an
$(n_3,n_3)$ configuration.
\ec

\proof
Indeed, in this case the lines are of the form $\star(v)$, for any $v
\in V$.
\qed

\be
\label{petdes}
If $G$ is the Petersen graph, then  {\rm Geo}$\, G$ is the Desargues
configuration. 
\ee

Indeed, let $G = (V,E)$ denote the Petersen graph, where the vertices
are described as the 2-subsets (written in the form $ij$) of $\hat{5}$
and $ij \edge kl$ is an edge if and only if $\{ i,j \} \cap \{ k,l\} =
\emptyset$:
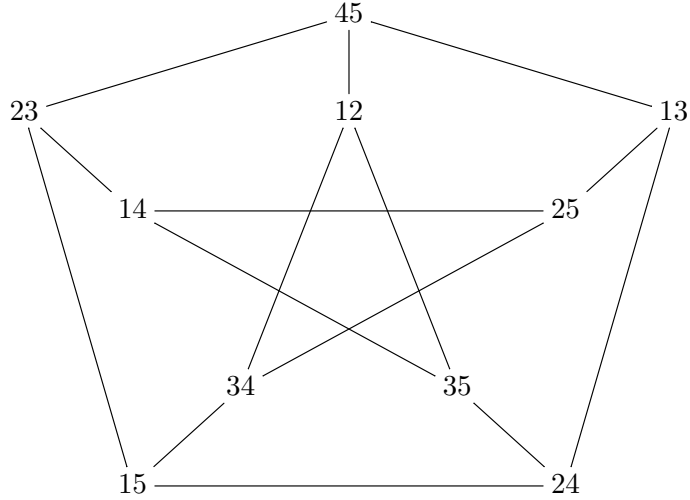
$$\xymatrix{
&&& 45 \ar@{-}[dlll] \ar@{-}[d] \ar@{-}[drrr] &&& \\
23 \ar@{-}[ddddr] \ar@{-}[dr] &&& 12 \ar@{-}[dddl] \ar@{-}[dddr]
&&& 13 \ar@{-}[ddddl] \ar@{-}[dl] \\
& 14 \ar@{-}[rrrr] \ar@{-}[ddrrr] &&&& 25 \ar@{-}[ddlll] & \\
&&&&&& \\
&& 34 \ar@{-}[dl] && 35 \ar@{-}[dr] && \\
& 15 \ar@{-}[rrrr] &&&& 24 &
}$$
Since the graph has girth 5, it follows easily that $\geo G = (V,\L)$
for $\L = \{ \star(v) \mid v \in V\}$, which coincides precisely with
our previous description of the Desargues configuration.

\bigskip

We say that a PEG $\G =
(P,\L)$ is {\em connected} if  there is no nontrivial partition $\L =
\L_1 \cup \L_2$ such that $(\cup \L_1)\cap(\cup \L_2) = \emptyset$. 
Note that this is equivalent to the usual geometric concept of
connectedness if our PEG 
has an euclidean geometric realization through real lines and real
points.

\bp
\label{cgeo}
Let $G$ be a graph in {\rm SC}$3$ with {\rm mindeg}$\, \G \geq
2$. Then the following conditions are equivalent:
\bi
\item[(i)] {\rm Geo}$\, G$ is connected;
\item[(ii)] $G$ is not bipartite.
\ei
\ep

\proof
By definition, $\geo G$ is disconnected if and only if there exists a
nontrivial partition $\L_G = \L_1 \cup \L_2$ such that $(\cup
\L_1)\cap(\cup \L_2) = \emptyset$. In view of Proposition \ref{peg}
and (G1), this supposes a nontrivial partition $V = V_1 \cup V_2$ with
$\cup \L_1 = V_1$ and $\cup \L_2 = V_2$. 

If $G$ is bipartite with respect to a partition $V = V_1 \cup V_2$,
then we take
$$\L_1 = \{ \star(v) \mid v \in V_2\}, \quad \L_2 = \{ \star(v) \mid v
\in V_1\}.$$ 
Since $\mindeg G \geq
2$, and by Proposition \ref{peg}, this shows that $\geo G$ is
disconnected. 

Conversely, assume that $\geo G$ is
disconnected. Hence there exists a nontrivial partition $V = V_1 \cup V_2$ with
$\cup \L_1 = V_1$ and $\cup \L_2 = V_2$. Suppose that $\star(v)
\subseteq V_1$ for some $v \in V_1$. Since $G$ is connected, it
follows easily from an induction argument that $\star(w)
\subseteq V_1$ for any $w \in V_1$, contradicting $V_2 \neq
\emptyset$. Hence $\star(v)
\subseteq V_2$ for every $v \in V_1$. By symmetry, we also have $\star(v)
\subseteq V_1$ for every $v \in V_2$. Therefore $G$ is bipartite.
\qed

As we mentioned in the beginning of the section, we can view graphs as
a particular case of PEG's, when we assume lines to have exactly two
points. 
Note that the concept of connectedness for PEG's coincides with the
usual concept of connectedness for graphs when we view graphs as PEG's.

Given a PEG $\G = (P,\L)$, we define the {\em
  Levi graph} of $\G$ \cite{Cox} by $\levi\G = (P \cup \L, E)$, where $E$
consists of the edges of the form $p \edge L$, for all $L \in \L$ and
$p \in L$. 

Viewing $K_3$ as a PEG, we have
$$
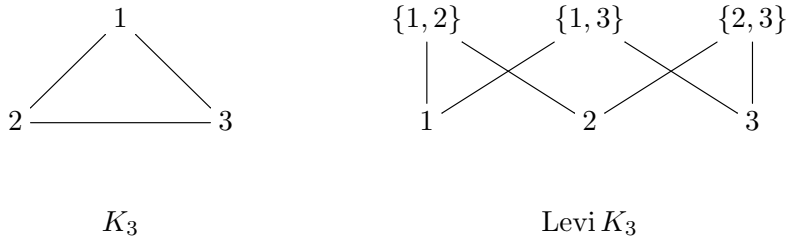
\xymatrix{
&1 \ar@{-}[dr] \ar@{-}[dl] &&&\{1,2\} \ar@{-}[d] \ar@{-}[dr] &\{1,3\}
\ar@{-}[dr] \ar@{-}[dl] &\{2,3\} \ar@{-}[dl] \ar@{-}[d]\\
2 \ar@{-}[rr] &&3&&1&2&3\\
&K_3&&&&\levi K_3&
}$$

If $G$ is a graph, its Levi graph is in fact a {\em subdivision} of
$G$. A simple way of 
picturing it is by introducing a new
vertex at the midpoint of every edge (breaking thus the original edge
into two). Obviously, the new vertices
represent the edges where they originated.

Among configurations, famous examples include the {\em Desargues
  graph} \cite{W5} as the Levi graph of 
the Desargues configuration  and the {\em Heawood graph} \cite{W1} as the Levi
graph of the Fano plane \cite{W7}.

The following results collects some elementary properties of the Levi
graph of a PEG (configuration) (see \cite{Cam,Cox}). Proofs
are immediate.

\bp
\label{levipeg}
Let $\G = (P,\L)$ be a PEG. Then:
\bi
\item[(i)] {\rm Levi}$\,\G$ is bipartite with respect to the partition $P
  \cup \L$;
\item[(ii)] the degree of $p \in P$ in {\rm Levi}$\,\G$ is the number
  of lines containing $p$;
\item[(iii)] the degree of $L \in \L$ in {\rm Levi}$\,\G$ is $|L|$;
\item[(iv)]  {\rm Levi}$\,\G$ has $|P|+|\L|$ vertices and
  $\displaystyle\frac{\sum_{p\in P}|\L(p)| \; + \; \sum_{L\in \L} |L|}{2}$ edges.
\ei
\ep

We define $\mindeg \G$ to be $\mindeg\levi G$.

\bc
\label{leviconf}
Let $\G = (P,\L)$ be an $(m_c,n_d)$ configuration. Then {\rm
  Levi}$\,\G$ has $m+n$ vertices and $cm = dn$ edges.
\ec

In particular, the Levi graph of the Desargues configuration, which is
a $(10_3,10_3)$ configuration, has 20 vertices and 30 edges.

Going back to the graph in Example \ref{coimbra}, it is easy to check
that $\geo G$ has $V = \hat{7}$ as set of points and lines $\star(v)$
for $v \in V$. The following picture shows that $\geo G$ is somehow
part of the Fano plane \cite{W7}:
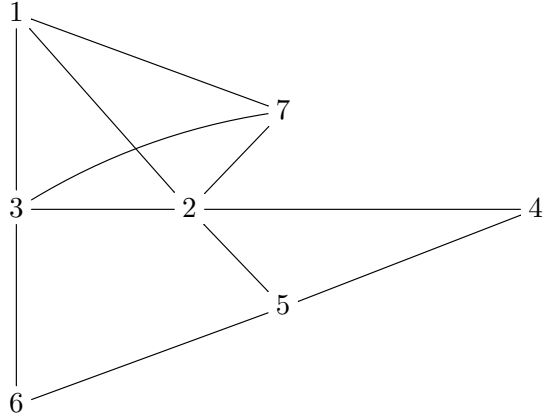
$$\xymatrix{
1 \ar@{-}[ddrr] \ar@{-}[dd] \ar@{-}[drrr] &&&&&& \\
&&& 7 \ar@{-}[dl] \ar@{-}@/_/[dlll] &&& \\
3 \ar@{-}[dd] \ar@{-}[rr] && 2 \ar@{-}[rrrr] \ar@{-}[dr] &&&& 4
\ar@{-}[dlll] \\
&&& 5 \ar@{-}[dlll] &&& \\
6 &&&&&&
}$$
Moreover, $\levi\geo \G$ can be obtained as follows: we make the Hasse
diagram of $\flats
G$ into a graph (the {\em Hasse graph} of $\flats G$) by taking as
vertices all flats, and letting $x \to y$ be an edge whenever $x$ covers $y$
in $\flats
G$ or vice-versa; removing the vertices $V$ and 
$\emptyset$, we get the {\em restricted Hasse graph} of $\flats G$,
which is then isomorphic to $\levi\geo \G$. This is just a particular
case of Proposition \ref{hasse}.

\medskip

We discuss next girth and connectedness.

\bp
\label{congirth}
Let $\G = (P,\L)$ be a PEG. Then
\bi
\item[(i)] {\rm gth}$\,${\rm Levi}$\,\G \geq 6$ and is not odd;
\item[(ii)] {\rm Levi}$\,\G$ is connected if and only if $\G$ is connected. 
\ei
\ep

\proof
(i) Since $\levi\G$ is bipartite by Proposition \ref{levipeg}(i), it has
no cycles of odd length. Therefore it is enough to exclude existence
of squares in $\levi\G$. Suppose that
$$
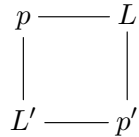
\xymatrix{
p \ar@{-}[r] \ar@{-}[d] & L \ar@{-}[d] \\
L' \ar@{-}[r] & p'
}$$
is a square in  $\levi\G$. Then $|L \cap L'| \geq 2$, contradicting
(G2). Therefore $\gth \levi\G \geq 6$.

(ii) Suppose that $\G$ is not connected. Then there is a
  nontrivial partition $\L = \L_1 \cup \L_2$ such that $(\cup
  \L_1)\cap(\cup \L_2) = \emptyset$. Suppose that $L \edge p \edge L'$
 is a path in $\levi\G$. Since $p \in L \cap L'$ and $(\cup
  \L_1)\cap(\cup \L_2) = \emptyset$, then $L$ and $L'$ must belong to
  the same side of the partition. Hence the connected component of a
  line in $\levi\G$ does not contain the lines in the other side of
  the partition, and so $\levi\G$ is not connected.

Conversely, suppose that $\levi\G$ is not connected. Let $\L_1$ be
the set of all lines in a fixed connected component of $\levi\G$ and let 
$\L_2 = \L \setminus \L_1$. Suppose that $p \in (\cup
  \L_1)\cap(\cup \L_2)$. Then there exist $L_1 \in \L_1$ and $L_2 \in
  \L_2$ such that $p \in L_1 \cap L_2$. Hence we have a path $L_1
  \edge p \edge L_2$ in $\levi\G$ and so $L_1$ and $L_2$ belong to
  the same connected component, a contradiction. Thus $(\cup
  \L_1)\cap(\cup \L_2) = \emptyset$ and so $\G$ is not connected.
\qed

Note that, if $G$ is a graph, the cycles of
$\levi G$ are of the form
$$v_0 \edge \{ v_0,v_1\} \edge v_1 \edge \{ v_1,v_2\} \edge \ldots
\edge \{ v_n,v_0\} \edge v_0,$$
whenever 
$$v_0 \edge v_1 \edge \ldots \edge v_n \edge v_0$$
is a cycle in $G$. Thus
$$\gth\levi G = 2\gth G.$$

\bp
\label{soberlevig}
The following conditions are equivalent for a PEG $\G = (P,\L)$:
\bi
\item[(i)] {\rm Levi}$\,\G$ is sober;
\item[(ii)] the mapping $P \to 2^{\L}: p \mapsto \L(p)$ is one-to-one;
\item[(iii)] for all distinct points $p,p' \in P$, there exists some
  line $L \in \L$ containing just one of them. 
\ei
\ep

\proof
We start by computing the stars of $\levi\G$. For $p \in P$ and $L
\in \L$, we have
$\star(p) = \L(p)$ and $\star(L) = L$ (recall that $L$ is a set of
points!). By axioms (G1) and (G3), we have repectively $\star(p) \neq
\emptyset$ and $\star(L) \neq
\emptyset$.
 Since $P \cap \L = \emptyset$, we must have always
$\star(p) \neq \star(L)$. On the other hand, the restriction
$\star|_{\L}$ is always one-to-one, hence $\levi\G$ is sober if and
only if $\star|_{\P}$ is one-to-one, which is equivalent to (ii). The
equivalence of (ii) and (iii) is trivial.
\qed

If $G$ is a graph, the above conditions are equivalent to saying that no
union of connected components of $G$ has exactly two vertices. 

We call a PEG satisfying the conditions of Proposition
\ref{soberlevig} {\em sober}. In view of axiom (G2), we immediately obtain:

\bc
\label{soberlevic}
If $\G$ is a PEG and {\rm mindeg}$\, \G \geq 2$, then $\G$ is sober. In
particular, if $\G$ is an $(m_c,n_d)$ configuration with $c \geq 2$,
then $\G$ is sober. 
\ec

This provides us with infinitely many examples of graphs in $\SC 3$
with girth $\geq 6$: 

\bc
\label{scpeg}
Let $\G$ be a PEG.
\bi
\item[(i)] If $\G$ is sober and connected, then {\rm Levi}$\,\G \in$
  {\rm\SC}$3$. 
\item[(ii)] If {\rm mindeg}$\, \G \geq 2$, then {\rm mindeg}$\,${\rm
    Levi}$\,\G \geq 2$. 
\ei
\ec

\proof
(i) Since $\G$ is sober, so is $\levi\G$. By Proposition \ref{congirth},
$\levi\G$ is connected and has girth $\geq 6$. Thus $\levi\G$ has
c-rank 3 by Corollary \ref{sclow}(iv).

(ii) By Proposition \ref{levipeg}, in view of $\mindeg \G \geq 2$ and (G3).
\qed

Note that, given a non bipartite cubic graph $C$ in $\SC 3$ with $n$
vertices (so $n \geq 10$), it follows from Proposition \ref{cgeo} and
Corollaries 
\ref{cuco} and \ref{soberlevic} that $\geo C$ is a sober connected $(n_3,n_3)$
configuration. Hence, by Proposition \ref{levipeg} and Corollary \ref{scpeg},
$\levi\geo C$
is now a bipartite 
cubic graph in $\SC 3$, so one can generate cubics this way. This does
not iterate because $\geo \levi\geo C$ does 
not stay connected. 

Given a graph $G = (V,E)$, we say that the vertex $v \in V$ is {\em
  closed} if $\{ v \} = \star(W)$ for some $W \subseteq V$, i.e.  $\{
v \} \in \flats G$. Note that this is also equivalent to the equality
$\{ v \} = \star(\star(v))$, since $\star(v)$ is clearly the greatest subset
$W$ of $V$ such that $v \in \star(W)$. We say that $G$ is {\em closed}
if all its vertices are closed.

By taking $G$ to be the graph $1 \edge 2 \edge 3 \edge 4$, and
omitting brackets/commas in the representation of the flats, we
can see that $\flats G = \{ 1234, 13, 24, 2, 3, \emptyset \}$ and so 2
and 3 are closed while 1 and 4 are not.

We can now prove the following (see \cite{CFS} in view of (ii)):

\bl
\label{uniatom}
Let $G = (V,E)$ be a finite graph satisfying one of the following two
conditions:
\bi
\item[(i)] $G$ is sober and cubic;
\item[(ii)] {\rm mindeg}$\, G \geq 2$ and $G$ has no squares.
\ei
Then $G$ is closed.
\el

\proof
Let $v \in V$. Clearly, $v \in \star(\star(v))$.  Suppose that $v \neq w \in
\star(\star(v))$. Then $\star(v) \subseteq \star(w)$. 

If $G$ is cubic, this implies $\star(v) = \star(w)$ and $G$ would not be
sober. Therefore (i) implies $\{ v \} =
\star(\star(v))$.

On the other hand, if (ii) holds, then by taking distinct $a,b \in
\star(v)$ we would get a square
$$\xymatrix{
v \ar@{-}[r] \ar@{-}[d] & a \ar@{-}[d] \\
b \ar@{-}[r] & w
}$$
a contradiction. Therefore we also get $\{ v \} =
\star(\star(v))$ in this case.
\qed

We can now prove the following result:

\bp
\label{hasse}
Let $G \in$ {\rm SC}$3$ have {\rm mindeg}$\, G \geq 2$. Then {\rm
  Levi}$\,${\rm Geo}$\, G$ is isomorphic to the restricted Hasse graph of {\rm
  Fl}$\, G$.
\ep

\proof 
Write $G = (V,E)$. By Proposition \ref{peg}, we have $\geo G = (V, \{
\star(v) \mid v \in V \})$ and so the vertex set of $\levi\geo \G$ is
$V \cup \{ \star(v) \mid v \in V \}$. On the other hand, by Lemmas
\ref{sct}(ii) and \ref{uniatom}(ii), the restricted Hasse graph $G'$ of
$\flats G$ has 
$$\{ \{ v\} \mid v \in V \} \cup \{ \star(v) \mid v \in V \}$$
as vertex set, yielding an obvious bijection to the vertex set of
$\levi\geo \G$. 

Now the edges of $\levi\geo \G$ are of the form $w \edge \star(v)$
whenever $w \in \star(v)$ $(v \in V)$, and this is equivalent to say
that $\star(v)$ covers $\{ w \}$ in $\flats G$. Therefore the two
graphs are isomorphic.
\qed

We proceed now to analyse the lattice of flats of the Levi graph of a
connected PEG with $\mindeg \geq 2$.

\bt
\label{flatslevig}
Let $\G = (P,\L)$ be a PEG and let {\rm Levi}$\,\G = (P \cup
\L,E)$. If $\G$ is connected and {\rm mindeg}$\, \G \geq 2$, then:
\bi
\item[(i)] {\rm Levi}$\,\G$ is closed;
\item[(ii)] {\rm Flats}$\,${\rm Levi}$\,\G = \{ P \cup \L, \emptyset\}
  \cup \{ \{ x \} \mid x \in P
  \cup \L \} \cup \{ L \mid L
  \in \L\} \cup \{ \L_p \mid p \in P\}$;
\item[(iii)] {\rm Flats}$\,${\rm Levi}$\,\G$ satisfies the
  Jordan-Dedekind condition.
\ei
\et

\proof
(i) By Lemma \ref{uniatom} and Proposition \ref{congirth}(i).

(ii) Given $p \in P$ and $L \in \L$, we have
$\star(p) = \L(p)$ and $\star(L) = L$. Moreover, $\star(p,L) =
\emptyset$. Now, given $p' \in P\setminus \{ p \}$, we have
$\star(p,p') = \{ L \}$ if $p,p' \in L \in \L$ (note that $L$ is then
unique by (G3)), otherwise empty. Finally, if $L' \in \L\setminus \{ L
\}$, we have in view of (G2) 
$\star(L,L') = \{ p \}$ if $L \cap L' = \{ p \}$, otherwise
empty. Note that we get all $\{ L \}$ by (G3) and (G2), and we  get
all $\{ p \}$ by (G1) and (G2). This proves (ii). 

(iii) Since $\mindeg \G \geq 2$, it follows easily from parts (i) and
(ii) that the maximal chains of $\flats \levi G$ are all of the form
$$\emptyset \subset \{ p \} \subset L \subset P\cup \L$$
or
$$\emptyset \subset \{ L \} \subset \L(p) \subset P\cup \L$$
for some $p \in L \in \L$. Therefore all maximal chains have
length 3.
\qed

We can now compute the c-independent subsets of $\levi\G$ for this
same class of PEG's:

\bc
\label{levindg}
Let $\G = (P,\L)$ be a PEG and let {\rm Levi}$\,\G = (P \cup
\L,E)$. If $\G$ is sober connected and {\rm mindeg}$\, \G \geq 2$, then
$W \subseteq P \cup \L$ is c-independent if and only if it satisfies one
of the following conditions:
\bi
\item[(i)] $|W| \leq 2$;
\item[(ii)] $|W| = 3$ and $|W \cap L| = 2$ for some $L \in \L$;
\item[(iii)] $|W| = 3$ and $|W \cap \L(p)| = 2$ for some $p \in P$.
\ei
\ec

\proof
By Theorems \ref{idind} and \ref{flatslevig}, $W$ is c-independent if
and only if $|W| \leq 2$ or
\beq
\label{levindg1}
|W| = 3,\; \star(W) = \emptyset \mbox{ and $W$ is not a potential line}.
\eeq
Thus we only need to show that the join of conditions (ii) and (iii)
is equivalent to (\ref{levindg1}).

Assume $|W| = 3$.
It is easy to see that $\star(W) \neq \emptyset$ can only occur if $W
\subseteq L$ for some $L \in \L$ or $W
\subseteq \L(p)$ for some $p \in P$. On the other hand, if $W$ is not a
potential line, then $|W \cap \star(x)| \geq 2$ for some $x \in
P\cup \L$, that is, either $|W \cap L| \geq 2$ or  $|W \cap \L(p)| \geq
2$ for some $p \in P$, $L \in \L$. Since $|W| = 3$, the
result follows.
\qed

Going back to the $K_3$ example at the beginning of this section, it
is now easy to check that every 3-subset $W$ of $V \cup E$ is
c-independent in $\levi K_3 \cong C_6$. Indeed, since $|E_v| = 2$ for
every $v \in V$, we only need to show that there exist necessarily
some $w_1,w_2 \in W$ at distance 2 (in $\levi K_3$). This is
certainly true for $C_6$, hence the c-independent subsets of vertices
of $\levi K_3$ (and therefore of $C_6$!) are all the subsets with at
most 3 vertices.

Another example is given by the Fano plane \cite{W7}. We have remarked before
that the Heawood graph $H$ is isomorphic to the  Levi
graph of the Fano plane and has no potential lines. It follows from
Theorem \ref{idind} that the c-independent subsets of $H = (V,E)$ are all
subsets with at most 3 vertices except the flats $\star(v)$ $(v \in
V)$. The reader can now check that these 463 subsets correspond to the
ones given by Corollary \ref{levindg}.

Given a PEG $\G = (P,\L)$, and since $\L \subseteq 2^P$, we can
consider the lattice $\widehat{\L}$ defined in Subsection
\ref{graphs}. We denote it by $\lat \G$. 

\bl
\label{isolat}
Given PEG's $\G$ and $\G'$ with $\mindeg \geq 2$, the following
conditions are equivalent: 
\bi
\item[(i)] $\G \cong \G'$;
\item[(ii)] {\rm Lat}$\, \G \cong$ {\rm Lat}$\, \G'$.
\ei
\el

\proof
It is immediate that the structure of $\G$ determines the structure of
$\lat \G$, up to isomorphism. Conversely, we can recover the structure
of $\G$ from $\lat \G$: 

Indeed, in view of (G2) and $\mindeg \G \geq 2$, we have
\beq
\label{isolat1}
\lat \G = \{ P, \emptyset\} \cup \L \cup P
\eeq
and so we can identify the points in $P$ with the atoms of $\lat \G$
and the lines in $\L$ with the maximal elements of $\lat\G \setminus \{
P\}$. Moreover $p \in L$ if and only if the corresponding atom of
$\lat \G$ is below the element representing in $L$, hence $\lat \G$
determines the structure of $\G$ up to isomorphism and the lemma
follows.
\qed  

If $\mindeg \G \geq 2$, we
can also introduce the {\em dual} PEG 
$\G^d$ (see \cite{Cox}):

\bl
\label{dualpeg}
Let $\G = (P,\L)$ be a PEG with {\rm mindeg}$\, \G \geq 2$. Then $\G^d = (\L,
\{ \L(p) \mid p \in P \})$ is also a PEG with {\rm mindeg}$\, \G^d \geq
2$. Moreover, {\rm Levi}$\,\G \cong$ {\rm Levi}$\,\G^d$.
\el

\proof
We have $\L \subseteq \cup_{p \in P} \L(p)$ since $\G$ satisfies
(G3). Hence  $\G^d$ satisfies
(G1). Given distinct $p,p' \in P$, we have $|\L(p) \cap \L(p')| \leq 1$ since
$\G$ satisfies (G2). Hence also $\G^d$ satisfies
(G2). Since $\mindeg \G \geq 2$ implies that $|\L(p)| \geq 2$ for
every $p \in P$, then $\G^d$ satisfies
(G3) and is thus a PEG.

Next, since $\G$ satisfies (G3), every $L \in \L$ belongs at least
to two $\L(p)$ and so $\mindeg \G^d \geq 2$.

Finally, let $\theta: P \cup \L \to \L \cup \{ \L(p) \mid p \in P \}$
be the bijection defined by $p\theta = \L(p)$ $(p \in P)$ and $L\theta
= L$ $(L \in \L)$. It is immediate that $\theta$ preserves the edges,
thus $\levi\G \cong \levi\G^d$.
\qed

Let $(X_1,\leq_1)$ and $(X_2,\leq_2)$ be lattices. We denote the maximum
(respectively the 
minimum) of both lattices by 1 (respectively 0) and assume the
remaining elements to be disjoint. The {\em coproduct} of $(X_1,\leq_1)$
and $(X_2,\leq_2)$, denoted by $(X_1,\leq_1) \sqcup (X_2,\leq_2)$, has
elements $X_1 \cup X_2$ (identifying the two 0's and the two 1's) and
partial order $\leq_1 \cup \leq_2$. In particular, $x_1 \wedge x_2 =
0$, $x_1 \vee x_2 = 1$ for all $x_1 \in X_1$ and $x_2 \in X_2$. 

\bt
\label{flpeg}
Let $\G$ be a PEG with {\rm mindeg}$\, \G \geq 2$. Then {\rm
  Fl}$\,${\rm Levi}$\,\G \cong$ {\rm Lat}$\,\G\; \sqcup$ {\rm
  Lat}$\,\G^d$. Moreover, this is the 
unique coproduct decomposition of  {\rm
  Fl}$\,${\rm Levi}$\,\G$.
\et

\proof
Write $\G = (P,\L)$. The isomorphism  $\flats
\levi\G \cong \lat\G \sqcup \lat\G^d$ follows easily from
Theorem \ref{flatslevig} and (\ref{isolat1}). 

Suppose now that $\p:\flats
\levi\G \to  X_1 \sqcup X_2$ is a lattice isomorphism for some
nontrivial lattices $X_1,X_2$. Let $Y_i$
denote the atoms of $\flats
\levi\G$ belonging to $X_i\p\inv$ $(i = 1,2)$. Suppose that $\{ p \}
\in Y_1$ with $p \in P$. If $p \edge L \edge p'$ is a path in
$\levi\G$, then it follows from (G2) that $\{ p \} \vee \{ p' \} =
L < P \cup \L$ and so $\{ p' \} \in Y_1$. Since $\G$ is connected, it
follows that $\{ q \} \in Y_1$ for every $q \in P$. Since $X_2$ is
nontrivial, then $\{ L \} \in Y_2$ for some $L \in \L$. If $L \edge q
\edge L'$ is a path in 
$\levi\G$, then it follows from (G2) that $\{ L \} \vee \{ L' \} \subseteq
\L(q) < P \cup \L$ and so $\{ L' \} \in Y_2$. Since $\G$ is connected, it
follows that $\{ M \} \in Y_2$ for every $M \in \L$. Since the atoms
determine the coproduct decomposition, it follows that $X_1 \cong
\lat\G$ and  $X_2 \cong
\lat\G^d$.
\qed

Now we can prove the following:

\bt
\label{isoflp}
Let $\G$ and $\G'$ be PEG's with {\rm mindeg}$\, \G$,
{\rm mindeg}$\, \G' \geq 2$. Then the following conditions are equivalent:
\bi
\item[(i)] $\G \cong \G'$ or $\G^d \cong \G'$;
\item[(ii)] {\rm Levi}$\,\G \cong$ {\rm Levi}$\,\G'$;
\item[(iii)]  {\rm
  Fl}$\,${\rm Levi}$\,\G \cong$ {\rm
  Fl}$\,${\rm Levi}$\,\G'$.
\ei
\et

\proof
(i) $\Rw$ (ii). In view of Lemma \ref{dualpeg}.

(ii) $\Rw$ (iii). Immediate.

(iii) $\Rw$ (i). Write $\G = (P,\L)$ and $\G' = (P',\L')$. Assume that
$\flats\levi\G \cong \flats\levi\G'$. By Theorem \ref{flpeg}, we
have $\lat(\G') \cong \lat(\G)$ or $\lat(\G') \cong \lat(\G^d)$. Now
(i) follows from Lemma \ref{isolat}.
\qed

The graph version is slightly simpler:

\bc
\label{isoflgra}
Let $G$ and $G'$ be finite connected graphs with  {\rm mindeg}$\, G$,
{\rm mindeg}$\, G' \geq 2$. Then the following conditions are equivalent:
\bi
\item[(i)] $G \cong G'$;
\item[(ii)] {\rm Levi}$\,G \cong$ {\rm Levi}$\,G'$;
\item[(iii)]  {\rm
  Fl}$\,${\rm Levi}$\,G \cong$ {\rm
  Fl}$\,${\rm Levi}$\,G'$. 
\ei
\ec

\proof
Viewing a graph $G$ as a PEG, its dual $G^d$ is a graph if and only if each
vertex of $G$ has degree 2, implying $G$ to be a cycle and therefore
self-dual. Now we apply Theorem \ref{isoflp}.
\qed

However, we recall that $\flats G \cong \flats G'$ does not imply $G
\cong G'$, even when $\mindeg G, \mindeg G' \geq 2$ (see the example
following the proof of Proposition \ref{soberrest}).

\section{Cubic graphs}

We present in this section some specific results concerning cubic graphs.

We start by some easy remarks concerning girth and c-rank of cubic graphs. 
For instance, note that $\gth G <
\infty$ for every finite cubic graph: any acyclic graph contains
necessarily vertices of degree 1.

In view of Propositions \ref{maxd} and \ref{trivials}, we have $2 \leq
\nk G \leq 4$ for every cubic graph $G$. However, the following result
shows that c-rank and girth are not independent for 
cubic graphs:

\bp
\label{cubic}
Let $G = (V,E)$ be a cubic graph.
\bi
\item[(i)] If ${\rm gth}\, G = 3$, then {\rm c}-{\rm rk}$\, G = 3$ or 4.
\item[(ii)] If ${\rm gth}\, G = 4$, then {\rm c}-{\rm rk}$\, G = 2$ or 3 or 4.
\item[(iii)] If ${\rm gth}\, G \geq 5$, then $G$ is sober,
{\rm c}-{\rm rk}$\, G = 3$ and $|\star(v,w)| \leq 1$ for distinct $v,w
\in V$.
\ei
Moreover, all these combinations with girth
$\leq 8$ can occur. If $G$ is sober and connected, only the cases
with ${\rm gth}\, G = 4$ and {\rm c}-{\rm rk}$\, G < 4$ are excluded.
\ep

\proof
(i) By Proposition \ref{casetwo}, since $G$ has a triangle.

(ii) By the comment preceding the proposition.

(iii) On the other hand,
Since ${\rm gth}\, G \geq 5$, $G$ has a restriction of the form
$$\bullet \edge \bullet \edge \bullet \edge \bullet$$
and so $\nk G > 2$ by Proposition \ref{casetwo}. On the other hand,
if $\nk G = 4$, then $G$ would have a square by Proposition
\ref{gfour}, a contradiction. Therefore $\nk G = 3$. Since $G$ has no
squares, the remaining conditions follow as well.

We present next examples to show that all these combinations with girth
$\leq 8$ occur:
\bi
\item
c-rank 2, girth 4: the complete bipartite graph $K_{3,3}$;
\item
c-rank 3, girth 3: 
$$
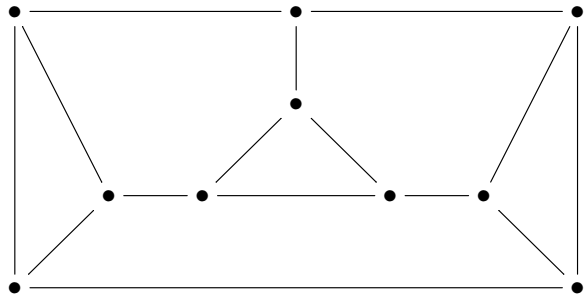
\xymatrix{
\bullet \ar@{-}[ddd] \ar@{-}[ddr] \ar@{-}[rrr] &&& \bullet \ar@{-}[rrr]
\ar@{-}[d] &&& \bullet \ar@{-}[ddl] \ar@{-}[ddd] \\
&&& \bullet \ar@{-}[dr] \ar@{-}[dl] &&& \\
& \bullet \ar@{-}[dl] \ar@{-}[r] & \bullet \ar@{-}[rr] && \bullet
\ar@{-}[r] & \bullet \ar@{-}[dr] & \\
\bullet \ar@{-}[rrrrrr] &&&&&& \bullet
}$$
\item
c-rank 3, girth 4:
$$
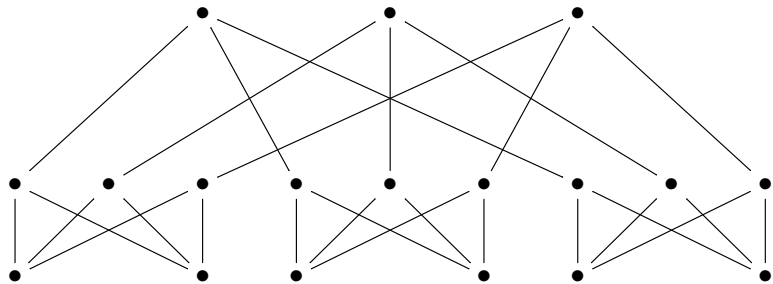
\xymatrix{
&&  \bullet \ar@{-}[ddll] \ar@{-}[ddr] \ar@{-}[ddrrrr] && \bullet
\ar@{-}[ddlll] \ar@{-}[dd] \ar@{-}[ddrrr] &&  \bullet \ar@{-}[ddllll]
\ar@{-}[ddl] \ar@{-}[ddrr] && \\
&&&&&&&& \\
\bullet \ar@{-}[d] \ar@{-}[drr] & \bullet \ar@{-}[dl] \ar@{-}[dr] &
\bullet \ar@{-}[dll] \ar@{-}[d] & \bullet
\ar@{-}[d] \ar@{-}[drr] & \bullet \ar@{-}[dl] \ar@{-}[dr] & \bullet
\ar@{-}[dll] \ar@{-}[d] & \bullet \ar@{-}[d]
\ar@{-}[drr] &
\bullet \ar@{-}[dl] \ar@{-}[dr] &  \bullet \ar@{-}[dll] \ar@{-}[d] \\
\bullet && \bullet & \bullet && \bullet & \bullet && \bullet
}$$
\item
c-rank 3, girth 5: the Petersen graph
\item
c-rank 3, girth 6: the Heawood graph \cite{W1};
\item
c-rank 3, girth 7: the McGee graph \cite{W2};
\item
c-rank 3, girth 8: the Tutte-Coxeter graph \cite{W3};
\item
c-rank 4, girth 3: the complete graph $K_4$;
\item
c-rank 4, girth 4:
$$
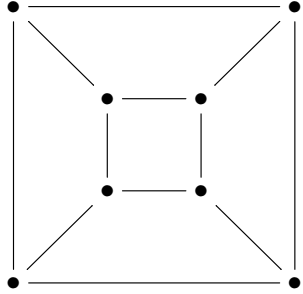
\xymatrix{
\bullet \ar@{-}[ddd] \ar@{-}[rrr]  \ar@{-}[dr] &&& \bullet \ar@{-}[ddd]
\ar@{-}[dl] \\ 
& \bullet \ar@{-}[d]  \ar@{-}[r] & \bullet \ar@{-}[d] & \\
& \bullet \ar@{-}[dl]  \ar@{-}[r] & \bullet \ar@{-}[dr] & \\
\bullet &&& \bullet  \ar@{-}[lll] 
}$$
\ei
Note that all these examples are sober and connected except those
with ${\rm gth}\, G = 4$ and $\nk G < 4$. The reason
for the exclusion of these combinations lies within Corollary
\ref{sclow}: if $G$ is sober and $\nk G < 4$, then $G$ has no squares.
\qed

Note that some of the arguments used in this proof are valid also for
graphs which are not cubic. For instance, if $\gth G \geq 5$ and all
vertices of $G$ have degree
$> 1$, then $G$ is necessarily sober.

It is an interesting problem to determine under which conditions the
lattice of flats of a graph has certain properties.

In the following theorems, we present results for the case of
connected cubic graphs. We start with a couple of useful lemmas.

\bl
\label{atoms}
Let $G = (V,E)$ be a finite nonempty graph. Then:
\bi
\item[(i)] every atom of {\rm Fl}$\, G$ is of the form {\rm St}$(${\rm St}$(v))$
  for some $v \in V$; 
\item[(ii)] the converse is true if $G$ is cubic.
\ei
\el

\proof
(i) Let $W$ be an atom of $\flats G$.
We may write $W = \star(X)$ for some $X \subseteq
V$. Let $v \in W  =
\star(X)$. Then $X \subseteq \star(v)$ and so $\star(\star(v))
\subseteq \star(X) = W$. Since $v \in \star(\star(v))$ and $W$ is an
atom, we get $\star(\star(v)) = W$.  

(ii) Assume that $G$ is cubic and $W = \star(\star(v))$ for some $v \in
V$. Let $u \in V$ be such that $W \cap \star(u) \neq \emptyset$. We
must prove that $W \subseteq \star(u)$. 

Indeed, if $x \in W \cap \star(u) =  \star(\star(v)) \cap \star(u)$,
then $\star(v) \cup \{ u \} \subseteq \star(x)$. Since $|\star(v)| =
|\star(x)| = 3$, it follows that $u \in \star(v)$ and so
$W = \star(\star(v)) \subseteq \star(u)$ as required. 
\qed

\bl
\label{3atom}
Let $G = (V,E)$ be a finite connected cubic graph. Then the following
conditions are equivalent:
\bi
\item[(i)] {\rm St}$(v)$ is an atom of {\rm Fl}$\, G$ for some $v \in V$;
\item[(ii)] $G \cong K_{3,3}$.
\ei
\el

\proof
(i) $\Rw$ (ii). 
If $\star(v)$ is an atom  of $\flats G$, then, for every $u \in V$, either
$\star(v) \subseteq 
\star(u)$ or $\star(v) \cap \star(u) = \emptyset$. Since $G$ is cubic,
$\star(v) \subseteq \star(u)$ is actually equivalent to $\star(v) =
\star(u)$. Writing $\star(v) = \{ a, b, c \}$ and $\star(a) = \{
v,x,y\}$, we can take $u$ above equal to $x$ and $y$ to obtain
$\star(x) = \star(y) = \star(v) = \{ a,b,c \}$. It follows that $G$
has a subgraph of the form
$$
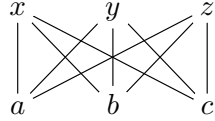
\xymatrix{
x \ar@{-}[d] \ar@{-}[dr]  \ar@{-}[drr] & y \ar@{-}[dl]
\ar@{-}[d] \ar@{-}[dr] & z \ar@{-}[d] \ar@{-}[dl]  \ar@{-}[dll] \\ 
a & b & c
}$$
(note that $\{ x,y,z \} \cap \{ a,b,c \} = \emptyset$ due to the
absence of loops). Since $G$ is cubic and connected, this must be the
whole of $G$, which is then isomorphic to $K_{3,3}$. 

(ii) $\Rw$ (i). Since the lattice of flats of $K_{3,3}$ is isomorphic
to $(2^{\hat{2}},\subseteq)$.
\qed

\bt
\label{latcc}
Let $G = (V,E)$ be a finite connected cubic graph. Then the following
conditions are equivalent:
\bi
\item[(i)]
{\rm Fl}$\, G$ is distributive;
\item[(ii)]
{\rm Fl}$\, G$ is modular;
\item[(iii)]
{\rm Fl}$\, G$ is semimodular;
\item[(iv)]
{\rm Fl}$\, G$ is geometric;
\item[(v)]
$G \cong K_4$ or $G \cong K_{3,3}$.
\ei
\et

\proof
The implications (i) $\Rw$ (ii), (ii) $\Rw$ (iii) and (iv) $\Rw$ (iii)
are immediate. Since the lattices of flats of $K_4$ and $K_{3,3}$ are isomorphic
respectively to $(2^{\hat{4}},\subseteq)$ and
$(2^{\hat{2}},\subseteq)$, we get (v) $\Rw$ (i) and (v) $\Rw$ (iv). It
remains to be proved that (iii) $\Rw$ (v).

Assume that $\flats G$ is semimodular. Suppose that $\diam G > 2$. Let
$v,w \in V$ be such that $d(v,w) > 2$. Write $\star(w) = \{ w_1,w_2,w_3\}$.
By Lemma \ref{atoms}, $\star(\star(w_j))$ is an atom of $\flats G$
for $j = 1,2,3$. Since $\diam K_{3,3} = 2$, it follows from Lemma
\ref{3atom} that $\star(v)$ is
not an atom of $\flats G$.
Write $\star(v) = \{ v_1,v_2,v_3\}$.
By Lemma \ref{atoms}, $\star(\star(v_i))$ is an atom of $\flats G$
for $i = 1,2,3$.

Suppose that
\beq
\label{latcc1}
\forall i,j \in \hat{3}\; \exists z_{ij} \in V: v_i,w_j \in
\star(z_{ij}).
\eeq
Since $d(v,w) > 2$, we must have $z_{ij} \neq v,w$. Hence $\{
z_{i1},z_{i2},z_{i3} \} \subseteq \star(v_i) \setminus \{ v \}$ and so
\beq
\label{latcc2} 
|\{ z_{i1},z_{i2},z_{i3} \}| \leq 2 \; \mbox{ for }i = 1,2,3.
\eeq
Similarly, 
\beq
\label{latcc3}
|\{ z_{1j},z_{2j},z_{3j} \}| \leq 2\; \mbox{ for }j = 1,2,3.
\eeq 

Let $G' = (V',E')$ be the graph such that $V' = \hat{3} \times
\hat{3}$ and $(i,j) \edge (i',j')$ is an edge if and only if  $(i,j)
\neq (i',j')$ and $z_{ij} = z_{i'j'}$. By (\ref{latcc2}) and
(\ref{latcc3}), $G'$ has at least 6 edges. Since $G'$ has 9 vertices,
there must be a pair of incident edges. Hence there exist distinct
$(i,j), (i',j'), (i'',j'') \in V'$ such that $z_{ij} = z_{i'j'} =
z_{i''j''}$. Thus  $v_i,v_{i'},v_{i''}, w_j,w_{j'},w_{j''} \in
\star(z_{ij})$ and so
$$|\{ v_i,v_{i'},v_{i''}, w_j,w_{j'},w_{j''} \}| \leq 3.$$
Since $\{ v_i,v_{i'},v_{i''} \} \cap \{ w_j,w_{j'},w_{j''} \} =
\emptyset$ due to $d(v,w) > 2$, this contradicts $(i,j), (i',j'),
(i'',j'')$ being all distinct. Therefore (\ref{latcc1}) fails and so
there exist $i,j \in \hat{3}$ such that $\star(v_i) \cap \star(w_j) =
\emptyset$. 

Now it is easy to check that 
$$
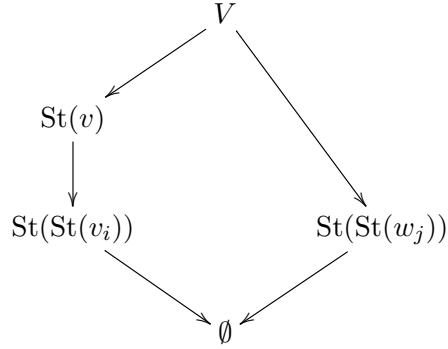
\xymatrix{
& V \ar[dl] \ar[ddr] & \\
\star(v) \ar[d] && \\
\star(\star(v_i)) \ar[dr] && \star(\star(w_j)) \ar[dl] \\
& \emptyset &
}$$
is a sublattice of $\flats G$. On the one hand, we have $\star(v) \cap
\star(\star(w_j)) \subseteq \star(v) \cap
\star(w) =
\emptyset$ since $d(v,w) > 2$. On the other hand, suppose that
$\star(\star(v_i)) \cup \star(\star(w_j)) \subseteq \star(z)$ for some $z \in
V$. Then $v_i,w_j \in \star(z)$, contradicting $\star(v_i) \cap \star(w_j) =
\emptyset$. This proves that $G$ cannot be semimodular if $\diam G >
2$. Hence $\diam G \leq 2$.

Suppose first that $G$ is not sober. Then we have $\star(a) = \{
x,y,z \} = \star(b)$ for some distinct $a,b \in V$, and so $G$ has a
subgraph of the form
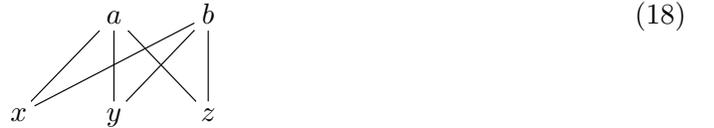
\beq
\label{new1}
\xymatrix{
& a \ar@{-}[dl] \ar@{-}[d] \ar@{-}[dr] & b \ar@{-}[dl] \ar@{-}[d]
\ar@{-}[dll] \\
x & y & z
}
\eeq
Supose that there exists an edge connecting two of the vertices
$x,y,z$, say $x \edge y$. Then $G$ has a subgraph of then form
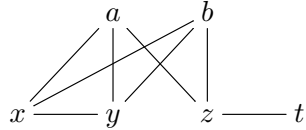
$$\xymatrix{
& a \ar@{-}[dl] \ar@{-}[d] \ar@{-}[dr] & b \ar@{-}[dl] \ar@{-}[d]
\ar@{-}[dll] & \\
x \ar@{-}[r] & y & z \ar@{-}[r] & t
}$$
and it is now clear that $d(t,y) > 2$, contradicting $\diam G \leq
2$. Hence (\ref{new1}) is a restriction of $G$. If there exists some
$X \in \flats G$ satisfying $ \{
x,y,z \} \supset X \supset \emptyset$, it is easy to check that
$$
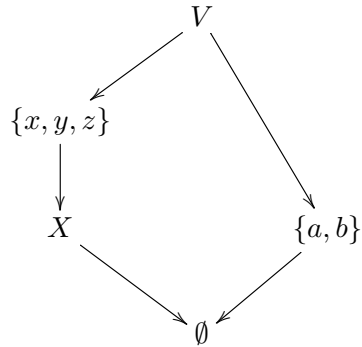
\xymatrix{
& V \ar[dl] \ar[ddr] & \\
 \{ x,y,z \} \ar[d] && \\
X \ar[dr] && \{ a,b \} \ar[dl] \\
& \emptyset &
}$$
would be a sublattice of $\flats G$, contradicting 
semimodularity. Thus $\{ x,y,z \}$ is an atom and so $G \cong
K_{3,3}$ by Lemma \ref{3atom}.

Therefore we may assume that $G$ is sober. Suppose first that there
exists some edge $a \edge b$ which does not lie in any triangle of
$G$. Since $G$ is sober, if follows from Lemma \ref{uniatom}(i) that
$\{ a \} = \star(\star(a))$ and $\{ b \}  =
\star(\star(b))$. Moreover, $a \in \star(v)$ for some $v \in V$. We claim
that 
$$\xymatrix{
& V \ar[dl] \ar[ddr] & \\
\star(v) \ar[d] && \\
\{ a \} \ar[dr] && \{ b \} \ar[dl] \\
& \emptyset &
}$$
is a sublattice of $\flats G$, a contradiction. 
Indeed, if $b \in \star(v)$, then $a,b,v$ would be the vertices of a
triangle, contradicting our assumption, and no flat can contain $a,b$
simultaneously by the same reason.
Hence every edge of
$G$ must lie in some triangle.

Now $G$ must have a subgraph of the form
$$
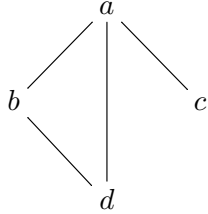
\xymatrix{
& a \ar@{-}[dl] \ar@{-}[dd] \ar@{-}[dr] & \\
b \ar@{-}[dr] && c \\
& d & 
}$$
Since the edge $a \edge c$ must lie in some triangle, we have an edge
$c \edge d$ or an edge $b \edge c$. Without loss of generality, we may
assume that $c \edge d$ is an edge.
If we have an edge $c \edge e$ with $e \neq b$, then we
have
$$
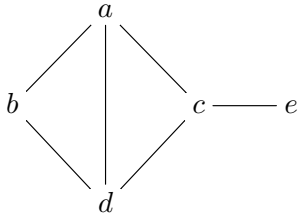
\xymatrix{
& a \ar@{-}[dl] \ar@{-}[dd] \ar@{-}[dr] && \\
b \ar@{-}[dr] && c \ar@{-}[dl] \ar@{-}[r] & e \\
& d &&
}$$ 
and so $d(b,e) = 3$, a contradiction. Hence there is an edge $b \edge
c$ as well and so $G \cong K_4$. Therefore (v) holds.
\qed

For every $n \geq 3$, we define the cylindrical strip $H_n$ by
$$
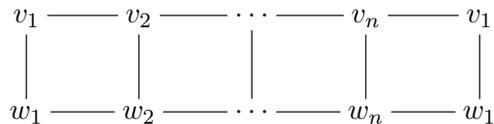
\xymatrix{
v_1 \ar@{-}[d] \ar@{-}[r] & v_2 \ar@{-}[d] \ar@{-}[r] & \cdots
\ar@{-}[d] \ar@{-}[r] & v_n \ar@{-}[r] \ar@{-}[d] & v_{1} \ar@{-}[d] \\ 
w_1 \ar@{-}[r] & w_2 \ar@{-}[r] & \cdots \ar@{-}[r] & w_{n} \ar@{-}[r] & w_1 
}$$
and the M\"obius strip $\tilde{H}_n$ by 
$$
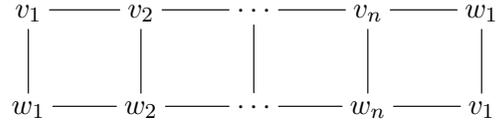
\xymatrix{
v_1 \ar@{-}[d] \ar@{-}[r] & v_2 \ar@{-}[d] \ar@{-}[r] & \cdots
\ar@{-}[d] \ar@{-}[r] & v_n \ar@{-}[r] \ar@{-}[d] & w_{1} \ar@{-}[d] \\ 
w_1 \ar@{-}[r] & w_2 \ar@{-}[r] & \cdots \ar@{-}[r] & w_{n} \ar@{-}[r] & v_1 
}$$

\bt
\label{jd}
Let $G = (V,E)$ be a finite connected cubic graph. Then the following
conditions are equivalent:
\bi
\item[(i)]
{\rm Fl}$\, G$ satisfies the Jordan-Dedekind condition;
\item[(ii)]
{\rm c}-{\rm rk}$\, G \leq 3$ or ($G$ is sober and every edge of $G$
lies in some square); 
\item[(iii)]
{\rm c}-{\rm rk}$\, G \leq 3$ or $G \cong K_4$ or $G \cong H_n$ for
some $n \geq 3$ 
or $G \cong \tilde{H}_n$ for some $n \geq 4$.
\ei
\et

\proof
(i) $\Rw$ (ii).
We may assume that $\nk G \geq 4$. Suppose that $G$ is not sober. Then
there exist $v,w \in V$ such that $\star(v) = \star(w)$. By Lemma
\ref{atoms}, $\star(\star(v))$ is an atom of $\flats G$. Since $v,w
  \in \star(\star(v))$, it follows that  $\flats G$ has an atom with
    2 elements. Since any $X \in \flats G \setminus \{ V \}$ has at
      most 3 elements, then $\flats G$ has a maximal chain with
        length $\leq 3$. Since $\nk G \geq 4$ implies the existence of
        some maximal chain with
        length 4, $G$ fails the Jordan-Dedekind condition.

Hence we may assume also that $G$ is sober. Suppose now that $a
\edge b$ is an edge of $G$. Write $\star(a) =
\{ b, c,d \}$. If $\star(b,c)$ contains some other element $x \neq a$,
then $a \edge b$ belongs to the square
$$
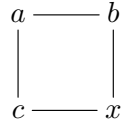
\xymatrix{
a \ar@{-}[r] \ar@{-}[d] & b \ar@{-}[d] \\
c \ar@{-}[r] & x
}$$
hence we may assume that $\star(b,c) = \{ a \}$ and so $\{ a\}$ is an atom of
$\flats G$. Write $\star(b) = \{ a,y,z \}$. If (i) holds, and since
$\nk G \geq 4$, the chain $\emptyset \subset \{ a \} \subset \{ a,y,z \}
\subset V$ must admit a refinement. We may therefore assume that $\{
a,y \} \in \flats G$. It follows that $\{
a,y \} = \star(p,q)$ for some distinct $p,q \in V$. Hence $p,q \in
\star(a) = \{ b, c,d \}$ and so $\{ p,q \} \cap \{ c,d \} \neq
\emptyset$. Assuming that $p \in \{ c,d \}$, we obtain a 4-cycle
$$
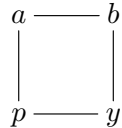
\xymatrix{
a \ar@{-}[r] \ar@{-}[d] & b \ar@{-}[d] \\
p \ar@{-}[r] & y
}$$
and so (ii) holds.

(ii) $\Rw$ (iii). We may assume that $\nk G \geq 4$, $G$ is sober and
every edge of $G$ lies in some square. 
We consider two cases:

\smallskip

\noindent
\underline{Case I}: $\gth G = 3$.

\smallskip

\noindent
Suppose first that $G$ has a subgraph $G_0$ the form 
$$
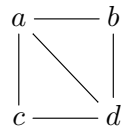
\xymatrix{
a \ar@{-}[r]  \ar@{-}[dr] \ar@{-}[d] & b \ar@{-}[d] \\
c \ar@{-}[r] & d
}$$
Then the edge $a \edge d$ must be part of a square
$$
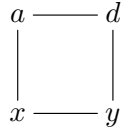
\xymatrix{
a \ar@{-}[r]  \ar@{-}[d] & d \ar@{-}[d] \\
x \ar@{-}[r] & y
}$$
Since $\star(a)$ and $\star(d)$ are fully determined, we must have $\{
x,y \} = \{ b,c\}$ and so $G$ has a subgraph isomorphic to
$K_4$. Since $G$ is connected and cubic, then $G \cong K_4$.

Hence we may assume that $G$ has no subgraph isomorphic to $G_0$
above. Take a triangle in $G$. Since every edge must belong to a
square and we are excuding subgraphs isomorphic to $G_0$, then $G$
must have a subgraph of the form
$$
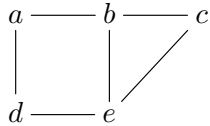
\xymatrix{
a \ar@{-}[r]  \ar@{-}[d] & b \ar@{-}[d] \ar@{-}[r] & c \ar@{-}[dl] \\
d \ar@{-}[r] & e & 
}$$
The existence of an edge $c \edge a$ or $c \edge d$ would imply the
presence of a 
subgraph isomorphic to $G_0$, hence we have an edge $c \edge f$ for
some new vertex $f$. Considering a square
$$
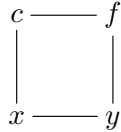
\xymatrix{
c \ar@{-}[r]  \ar@{-}[d] & f \ar@{-}[d] \\
x \ar@{-}[r] & y 
}$$
it follows easily that either $x = b$ and $y = a$, or $x = e$ and $y =
d$. These cases yield in fact isomorphic subgraphs, hence we assume
the first to get
$$
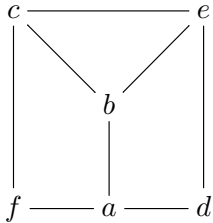
\xymatrix{
c \ar@{-}[rr] \ar@{-}[dr] \ar@{-}[dd] && e \ar@{-}[dd]
\ar@{-}[dl] \\
& b \ar@{-}[d] &  \\
f \ar@{-}[r] & a  \ar@{-}[r] & d
}$$
It is straightforward to check that the only square that can contain
the edge $c \edge e$ is 
$$
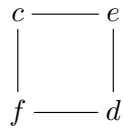
\xymatrix{
c \ar@{-}[r]  \ar@{-}[d] & e \ar@{-}[d] \\
f \ar@{-}[r] & d
}$$
hence $G$ contains a subgraph isomorphic to $H_3$ and is therefore
isomorphic to $H_3$. 

\smallskip

\noindent
\underline{Case II}: $\gth G = 4$.

\smallskip

\noindent
Let $G'$ be a subgraph of $G$ of the form
$$
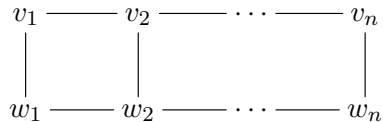
\xymatrix{
v_1 \ar@{-}[r] \ar@{-}[d] & v_2 \ar@{-}[r] \ar@{-}[d] & \cdots
\ar@{-}[r] & v_n  \ar@{-}[d] \\
w_1 \ar@{-}[r] & w_2 \ar@{-}[r] & \cdots
\ar@{-}[r] & w_n
}$$
with $n$ maximum. We claim that $n \geq 4$.

Indeed, suppose first that $n = 2$. Since $G$ has no triangles, then we have
a subgraph
$$
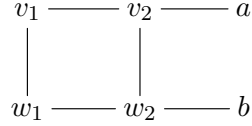
\xymatrix{
v_1 \ar@{-}[r] \ar@{-}[d] & v_2 \ar@{-}[r] \ar@{-}[d] & a \\
w_1 \ar@{-}[r] & w_2 \ar@{-}[r] & b 
}$$
Let 
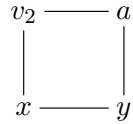
$$\xymatrix{
v_2 \ar@{-}[r] \ar@{-}[d] & a \ar@{-}[d] \\
x \ar@{-}[r] & y
}$$
be a square containing $v_2 \edge a$.
Then either $x = v_1$ or $x = w_2$. Suppose first that $x = v_1$. If
$y = w_1$, then $\star(v_2) = \star(w_1)$, contradicting $G$ being
sober. On the other hand, if $y$ is a new vertex $c$, we get two
adjacent squares and contradict the maximality of $n$. The case $x =
w_2$ is similar, hence $n > 2$ in this case.

Suppose now that $n = 3$. Since there are no triangles and $G$ is
sober, we have a subgraph
$$
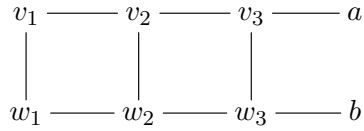
\xymatrix{
v_1 \ar@{-}[r] \ar@{-}[d] & v_2 \ar@{-}[r] \ar@{-}[d] & v_3 \ar@{-}[r]
\ar@{-}[d] & a \\
w_1 \ar@{-}[r] & w_2 \ar@{-}[r] & w_3 \ar@{-}[r] & b 
}$$
To avoid contradicting the maximality of $n$, we cannot accept an edge
$a \edge b$. Considering squares containing the edges $v_3 \edge a$
and $w_3 \edge b$, we obtain edges $v_1 \edge a$ and $w_1 \edge
b$. Taking an edge $a \edge c$, where $c$ is necessarily a new vertex,
we immediately get a contradiction by trying to fit the new edge into
a square. Thus $n \geq 4$. 

Now if $v_n \edge a$ is an edge, where $a$ is a new vertex, we cannot
fit this edge into a square without compromising the maximality of $n$,
hence we have either edges $v_1 \edge v_n$ and $w_1 \edge w_n$
(yielding $H_n$) or edges $v_1 \edge w_n$ and $w_1 \edge v_n$
(yielding $\tilde{H}_n$). Therefore (iii) holds.

(iii) $\Rw$ (ii). Immediate.

(ii) $\Rw$ (i). The case $\nk G = 2$ being trivial, suppose first
that $\nk G = 3$. Since the flats $\star(v)$ are the maximal elements of
$\flats G \setminus \{ V \}$ and no such flat is an atom of
$\flats G$ in view of 
Lemma \ref{3atom}, it follows that every maximal chain of $\flats G$
must have length 3 and so $G$ satisfies the Jordan-Dedekind condition.

Finally, assume that $\nk G > 3$. Let $v \in V$. By Lemma
\ref{uniatom}(i), $\{ v \}$ is an atom for every $v \in V$. Now, if $\{ a \}
\subset \{ a,b,c \} = \star(x)$ is 
a chain in $\flats G$, then we may assume that there exists some square 
$$
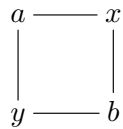
\xymatrix{
a \ar@{-}[r] \ar@{-}[d] & x \ar@{-}[d] \\
y \ar@{-}[r] & b
}$$
Since $G$ is sober, it follows that $\star(x,y) = \{ a,b\}$ and so
all maximal chains in $\flats G$ must have length 4. Thus (i) holds.
\qed

It is easy to check that all graphs $H_n$ and $\tilde{H}_n$ are
{\em vertex-transitive}: for all vertices $v$ and $w$, there exists an
automorphism $\p$ of the graph such that $v\p = w$ (i.e. all vertices
lie in a single automorphic orbit). 

By Proposition \ref{maxd} and Corollary \ref{sclow}, a finite sober
connected cubic graph has c-rank 4 if and only if it has a square. If
it is also vertex-transitive, then every vertex must lie in some square. 
The next example shows that one cannot replace {\em edge} by {\em
  vertex} in condition (ii) of Theorem \ref{jd}, even if we 
require vertex-transitivity:
$$
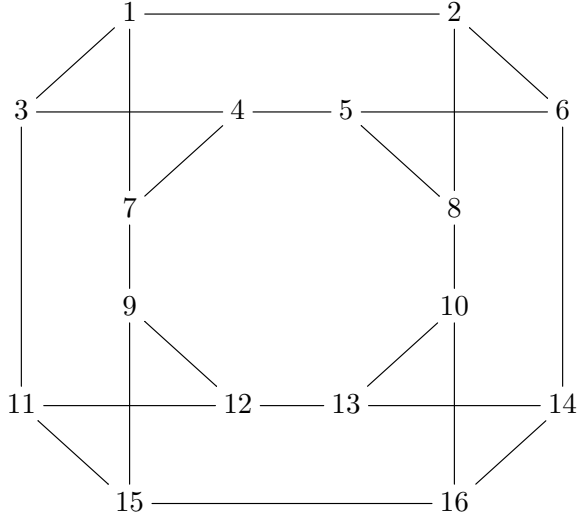
\xymatrix{
& 1 \ar@{-}[rrr] \ar@{-}[dl] \ar@{-}[dd] &&& 2 \ar@{-}[dr]
\ar@{-}[dd] & \\
3 \ar@{-}[rr] \ar@{-}[ddd] && 4 \ar@{-}[r] \ar@{-}[dl] &
5 \ar@{-}[rr] \ar@{-}[dr] && 6 \ar@{-}[ddd] \\
& 7 \ar@{-}[d] &&& 8 \ar@{-}[d] & \\
& 9 \ar@{-}[dd] \ar@{-}[dr] &&& 10 \ar@{-}[dd] \ar@{-}[dl] & \\
11 \ar@{-}[rr] \ar@{-}[dr] && 12 \ar@{-}[r] & 13
\ar@{-}[rr] && 14 \ar@{-}[dl] \\
& 15 \ar@{-}[rrr] &&& 16 & 
}$$
Indeed, this finite sober connected cubic graph has c-rank 4 and is
vertex-transitive (hence every vertex lies in some square), and yet it
fails the conditions of Theorem \ref{jd}. Note that in this case
$$V \supset \star(1) \supset \star(1,4) \supset \star(1,4,11) \supset
\emptyset$$
and
$$V \supset \star(3) \supset \star(3,12) \supset
\emptyset$$
are two maximal chains of different length.

\section{Minors and cm-rank}

Recall that a finite graph $G'$ is said to be a {\em minor} of a
finite 
graph $G$ if $G'$ can be obtained (up to isomorphism)
from $G$ by 
successive application of the following three operations:
\bi
\item[(D$_1$)] {\em vertex-deletion}: we delete a vertex;
\item[(D$_2$)] {\em edge-deletion}: we delete an edge;
\item[(C)] {\em contraction}: we delete an edge $v \edge w$ and
  identify the vertices $v$ and $w$.
\ei
If $G'$ is a minor of $G$, we write $G' \leq_m G$. 

It is easy to check that these operators commute with each other in
the sense that 
$$\begin{array}{c}
D_1D_2(G) \subseteq (D_2D_1 \cup D_1)(G),\quad
D_2D_1(G) \subseteq D_1D_2(G),\\
CD_1(G) \subseteq D_1C(G),\quad
D_1C(G) \subseteq (CD_1 \cup D_1^2)(G),\\
CD_2(G) \subseteq (D_2C \cup C)(G),\quad
D_2C(G) \subseteq (CD_2 \cup CD_2^2)(G),
\end{array}$$
hence a minor
of $G$ can in particular be obtained by applying to $G$ sequences of
contractions followed by edge-deletions followed by
vertex-deletions. Clearly, c-rank cannot increase by means of vertex-deletions
since we are bound to get a submatrix of he original one. However, the
example following  
Proposition \ref{res} shows that c-rank can increase
by means of edge-deletions. The same
happens for contractions: taking the very same square as an example,
which has c-rank 2, and performing a contraction, we get $K_3$ which
has higher c-rank.

Thus we introduce a second rank function for finite graphs: given a
finite connected graph $G$, let
$$\cmr G = \max\{ \nk G' \mid G' \leq_m G \}.$$
Since a minor has at most as many vertices as the original graph,
cm-rank is well defined. For every $m \in \N$, we denote by $\G_m$ the
class of all finite graphs with cm-rank $\leq m$. Since the minor
relation is transitive, $\G_m$ is closed for minors. In view of the
Robertson-Seymour Theorem (see \cite{Die}), there exists a finite set
of graphs $\F$ such that
$$G \in G_m \iff \forall F \in \F \hspace{.3cm} F \not\leq_m G.$$
We can easily construct the set $\F$ of forbidden graphs in our
case. For $m \geq 1$, let $\F_m$ consist of representatives of all
isomorphism classes of graphs with at most
$2m$ vertices and c-rank $m+1$. Let $\F_0$ contain a one-vertex graph.

\bp
\label{avoidRS}
The following conditions are equivalent for every finite graph $G$ and
every $m \in \N$:
\bi
\item[(i)] $G \in \G_m$;
\item[(ii)] $\forall F \in \F_m \; F \not\leq_m G$.
\ei
\ep

\proof
The case $m = 0$ holding trivially, we assume that $m \geq 1$.

(i) $\Rw$ (ii). 
If $G$ has a minor $G' \cong F \in F_m$, then $\cmr G \geq \nk G' =
\nk F = m+1$ and so $G \notin \G_m$.

(ii) $\Rw$ (i). If $G \notin \G_m$, then $G$ has a minor $G' = (V',E')$ of
c-rank $> m$. Since a subgraph of a minor is itself a minor, we may
assume that $\nk G' = m+1$. Hence there exist $I,J \subseteq V'$ such
that $|I| = |J| = m+1$ and $A_{G'}^c[I,J]$ is nonsingular. Write $I =
\{ i_1,\ldots,i_{m+1} \}$ and $J =
\{ j_1,\ldots,j_{m+1} \}$. Reordering
rows and columns if necessary, we may 
assume that $A^c[I,J]$ is of the form
(\ref{nons1}), for the ordering $i_1 < \ldots < i_{m+1}$ and $j_1 <
\ldots < j_{m+1}$. Replacing $j_1$ by $i_1$ and $i_{m+1}$ by
$j_{m+1}$, the resulting matrix is still of the form
(\ref{nons1}). Let $F$ be the restriction of $G'$ induced by the
vertices $\{ i_1,\ldots,i_{m}, j_2,\ldots,j_{m+1} \}$. Up to
isomorphism, we have $F \in \G_m$. Since $F \leq_m G' \leq_m G$, (ii)
fails as required.
\qed

Next we initiate a discussion on how the computation of cm-rank relates
to the matrix representation of graphs. A sequence of contractions
on a graph $G = (V,E)$ determines a partition $P:V = V_1 \cup
\ldots \cup V_m$
corresponding to the subsets of vertices that are eventually
identified into a single one. It is immediate that the restriction of
$G$ induced by each $V_i$ must be connected (we call such a partition
{\em connected}). How do we identify a
connected restriction within $A^c$? Through the following straightforward
observation:

\bp
\label{connmat}
The following conditions are equivalent for a finite graph $G = (V,E)$:
\bi
\item[(i)] $G$ is connected;
\item[(ii)] there exists no nontrivial partition $V = I \cup J$
  such that $A[I,J]$ is the null matrix;
\item[(iii)] there exists no nontrivial partition $V = I \cup J$
  such that all the entries in $A^c[I,J]$ are equal to 1.
\ei
\ep

What happens to the adjacency matrix when we perform a sequence of
contractions inducing the partition $P:V = V_1 \cup \ldots \cup V_m$?
Let the new graph be $G/P = (V/P,E/P)$, with $V/P = \hat{m}$, where
each vertex $i$
corresponds to the identification of the vertices in $V_i$. It is
straightforward to check that
$$A^c_{G/P}[i,j] = \left\{
\begin{array}{ll}
0& \mbox{ if $i \neq j$ and 0 occurs in }A^c_G[V_i,V_j]\\
1& \mbox{ otherwise}
\end{array}
\right.$$
If we follow a sequence of contractions by a sequence of
edge-deletions, we are entitled to replace 0s by 1s in the matrix
$A^c_{G/P}$. Finally, vertex-deletions correspond to deleting rows and
columns in this modified matrix, which does not increase c-rank, and
can therefore be ignored in the computation of the cm-rank. We
therefore obtain: 

\bp
\label{matcon}
Let $G = (V,E)$ be a finite graph. Then
{\rm cm}-{\rm rk}$\, G$ is the maximum value of {\rm rk}$\, \wt{A^c_{G/P}}$
when $P: V = V_1 \cup 
\ldots \cup V_m$ is a connected partition of $V$ and  
$$\wt{A^c_{G/P}}[i,j] = \left\{
\begin{array}{cl}
0\mbox{ or }1& \mbox{ if $i \neq j$ and 0 occurs in }A^c_G[V_i,V_j]\\
1& \mbox{ otherwise}
\end{array}
\right.$$
\ep

\section{The complement graph}

Given a graph $G = (V,E)$, its {\em complement graph} $\oo{G} = (V,\oo{G})$ is
the graph defined by the condition
$$\{ v,w \} \in \oo{E} \hspace{.5cm} \iff  \hspace{.5cm} \{ v,w \} \notin E\},$$
for all distinct $v,w \in V$. 

The classical idea of independence for a subset $W$ of vertices of $G$
(no edges between them) is related to our notion of c-independence by
$W$ being necessarily c-independent in $\oo{G}$, but not conversely.

We can get a lower bound for $\nk \oo{G}$ through the chromatic number.
An {\em edge coloring} of a graph $G = (V,E)$ with $c$ colors is a partition
$V = V_1 \cup \ldots \cup V_c$ such that no edge of $G$ connects two
vertices in the same $V_j$.
The {\em chromatic number} $c(G)$ is the minimum number $c$ of colors
to edge color $G$. 

\bp
\label{chro}
Let $G$ be a finite graph. Then {\rm c-rk}$\, \oo{G} \geq
\frac{|V|}{c(G)}$.
\ep

\proof
Since $|V_j| \geq \frac{|V|}{c(G)}$ for some $j$, then $\oo{G}$ has a
complete subgraph with at least $\frac{|V|}{c(G)}$ vertices and the
claim follows from Proposition \ref{res}(ii).
\qed

An important issue consists
of the study of the sum $S = \nk G + \nk \oo{G}$ for a graph with $n$
vertices. The examples we analyzed so far show that $S$ can be as
small as $\frac{n}{2} + 2$ (taking $G = K_{\frac{n}{2},\frac{n}{2}}$
for $n$ even, then $\oo{G}$ is a disjoint union of two copies of
$K_{\frac{n}{2}}$) and as large as $n+2$ (taking $G$ to be a graph
of the form 
$$v_1 \edge v_2 \edge \ldots \edge v_n$$
for $n \geq 4$). The next result offers an upper bound for $S$:

\bp
\label{boundsum}
Let $G = (V,E)$ be a finite graph with $|V| = n$. Then $S = $ {\rm
  c}-{\rm rk}$\, G \; +$ {\rm c}-{\rm rk}$\, \oo{G} < \sqrt{2} n +1$.
\ep

\proof
Assume that $|E| = k$ and $\nk G = m$. Then $m = \rk A^c_G$ and 
the witness characterization in Proposition \ref{altrank} yields 
$$(m-1)+\ldots +2+1 = \frac{m(m-1)}{2} \leq k.$$
Thus
$2k \geq m^2-m$, yielding $m \leq \frac{1+\sqrt{1+8k}}{2}$. 

Similarly, since $\oo{G}$ has $\frac{n(n-1)}{2}-k$ edges, we get 
$$\nk \oo{G} \leq \frac{1+\sqrt{1+4n^2-4n -8k}}{2}.$$
Hence 
$$S \leq 1+ \frac{\sqrt{1+8k}+ \sqrt{1+4n^2-4n -8k}}{2}.$$
A simple calculus exercise shows that a real-valued function of the
form
$$\sqrt{1+8x}+ \sqrt{1+8(a-x)} \quad (a > 0)$$
reaches its maximum when $x = \frac{a}{2}$. Hence
$$S \leq 1+ \sqrt{1+2n^2-2n} \leq 1+ \sqrt{2}n.$$
\qed

We can also note the following:

\bp
\label{ramsey}
Let $(G_n)_n$ be a sequence of nonisomorphic finite graphs and let $M_n =
\max\{$ {\rm c}-{\rm rk}$\, G_n,$ {\rm c}-{\rm rk}$\, \oo{G}_n \}$. Then 
$$\displaystyle \lim_{n\to +\infty} M_n = +\infty.$$
\ep

\proof
Let $R(k,k)$ denote the Ramsey number that ensures every complete
graph with at least $R(k,k)$ vertices, with edges colored by two
colors, to have a monocromatic complete subgraph with $k$
vertices. Let $k \in \mathbb{N}$. Since the graphs $G_n$ are
nonisomorphic, there exists some $p \in \N$ such that all graphs $G_n$
have at least $R(k,k)$ vertices for $n > p$. In particular, either
$G_n$ or $\oo{G}_n$ must contain a complete subgraph with $k$
vertices, and so $M_n \geq k$ by Proposition \ref{res}(ii). Therefore
$\lim_{n\to +\infty} M_n = +\infty.$
\qed 

We can give another perpective of the complement graph through the
{\em dual lattice of closed stars}. Given a graph $G = (V,E)$, the {\em
  closed star} of a vertex $v \in V$ is defined by
$$\cstar(v) = \star(v) \cup \{ v \}.$$

Given $\S \subseteq 2^V$,
it is easy to see that
$$\dwh{\S} = \{ \cup S \mid S \subseteq \S \}$$
is the $\vee$-subsemilattice of $(2^V, \subseteq)$ generated by
$\S$. Note that $\cup \S = \max \dwh{\S}$, and also $\emptyset =
\cup \emptyset = \min \dwh{\S}$. Similarly to the dual case,
$(\dwh{\S},\subseteq)$ is 
itself a lattice with
$$P \wedge Q = \cup \{ X \in \S \mid P\cap Q \subseteq X\}.$$
In particular, we can take 
$\oo{\S}_V = \{ \cstar(v) \mid v \in V
\}$ and consider the lattice $\dwh{\oo{\S}_V}$, which we call  the
{\em dual lattice of closed stars} of $G$.

\bt
\label{dualstars}
Let $G = (V,E)$ be a finite graph. Then {\rm c}-{\rm rk}$\, \oo{G} =$ {\rm
  ht}$\, \dwh{\oo{\S}_V}$.
\et

\proof
We know that $\nk \oo{G}$ is the maximum length $n$ of a chain of the form
\beq
\label{dualstars1}
V \supset \star_{\oo{G}}(v_1) \supset \star_{\oo{G}}(v_1,v_2)
\supset \ldots \supset \star_{\oo{G}}(v_1,\ldots,v_n) = \emptyset.
\eeq
Now
$$\star_{\oo{G}}(v_1,\ldots,v_i) = \star_{\oo{G}}(v_1) \cap \ldots
\cap \star_{\oo{G}}(v_i) = (V\setminus \cstar_{G}(v_1)) \cap \ldots
\cap (V\setminus \cstar_{G}(v_i))$$
and so 
$$V\setminus \star_{\oo{G}}(v_1,\ldots,v_i) = \cstar_{G}(v_1) \cup
\ldots \cup \cstar_{G}(v_1).$$
Passing (\ref{dualstars1}) to complement, it follows that  $\nk
\oo{G}$ is the maximum length $n$ of a chain of the form 
$$
\emptyset \subset \cstar_{G}(v_1) \subset \cstar_{G}(v_1) \cup
\cstar_{G}(v_2) \subset \ldots \subset \cstar_{G}(v_1) \cup \ldots \cup
\cstar_{G}(v_n) = V,
$$
which is precisely $\het \dwh{\oo{\S}_V}$.
\qed

As an example, we can now apply this result to the computation of the
c-rank of the complement of the Petersen graph $P$:

\be
\label{comppet}
{\rm c}-{\rm rk}$\, \oo{P} = 5$.
\ee

Write $P = (V,E)$, $\star(v) = \star_P(v)$ and $\cstar(v) =
\cstar_P(v)$. Assume that 
$$
\emptyset \subset \cstar(v_1) \subset \cstar(v_1) \cup
\cstar(v_2) \subset \ldots \subset \cstar(v_1) \cup \ldots \cup
\cstar(v_n) = V,
$$
is a chain of maximum length in $\dwh{\oo{\S}_V}$. We claim that
\beq
\label{comppet1}
|\cstar(v_1) \cup \cstar(v_2) \cup \cstar(v_3)| \geq 8.
\eeq

Since $P$ is
cubic, it follows from Proposition \ref{cubic}(iii) that $|\star(v_1)
\cup \star(v_2)| \geq 5$. 
If $v_1 \edge v_2$ is not an edge of $P$, then $|\cstar(v_1)
\cup \cstar(v_2)| \geq 7$ and so (\ref{comppet1}) must hold in this case.
Hence we may assume that $v_1$ and $v_2$ are adjacent in $P$.
Since $P$ has no triangles, then $\star(v_1,v_2) = \emptyset$ and so
$|\cstar(v_1) \cup \cstar(v_2)| = |\star(v_1) \cup \star(v_2)| = 6$.
Suppose that $|\cstar(v_1) \cup \cstar(v_2) \cup \cstar(v_3)| <
8$. Then $|\cstar(v_3) \cap \cstar(v_i)| \geq 2$ for some $i \in \{
1,2 \}$. If $v_i \edge v_3$ is an edge, we get a triangle in $P$; if
$v_i \edge v_3$ is not an edge, we get a square in $P$, a
contradiction in any case since $\gth P = 5$. Therefore
(\ref{comppet1}) holds and so $n \leq 5$.

It is a simple exercise to produce a chain of length 5 in
$\dwh{\oo{\S}_V}$, hence $\nk \oo{P} = 5$.

\section{Open questions}

Here is a list of open questions, some mentioned in the preceding
text:
\begin{enumerate}
\item
Characterize those finite graphs $G$ whose lattice of flats are
distributive, modular, semimodular, satisfy the Jordan-Dedekind chain
condition, or are extremal lattices. 

For extremal lattices see \cite{Mar}.
The most important questions are for the Jordan-Dedekind chain
condition and for the extremal lattices which need not satisfy the
Jordan-Dedekind chain condition. See Theorems \ref{latcc} and \ref{jd} for some
results. 
\item
Which hereditary collections have Boolean representations? 

See Subsection \ref{matroi} and \cite{IR1} for definitions.  All
matroids have boolean representations (see \cite{IR2}), but what else? 
This will be the topic of the future paper \cite{RSil}, which will
also carry out a similar analysis like Section \ref{gvx} but for
higher c-rank.   
\item
When do two graphs have isomorphic lattices of flats?
\item
Which matroids arise as the c-independent sets of a graph? Prove not
all matroids arise this way. 
\item
When is $\geo G$, for (cubic) $G \in \SC 3$ realizable as lines in
Euclidean space (i.e. realizable affinely)? 
Analyse  further the map which associates to a non bipartite cubic $C \in
\SC 3$ the bipartite cubic $\levi\geo C$.
\item
Compare the results of this paper with those of Brijder and Traldi in
\cite{BT}. 
\item
Extend the analysis of graphs in this paper to Moore graphs,
generalized Petersen graphs, etc. What is $\matro M$ when $M$ is the
Moore mystery graph of girth 5 (which may or may not exist) on 57
vertices? See \cite{W8}. 
\item
What is the smallest number of edges we can have in a graph with
$n$ vertices to maximize (minimize) $\nk G + \nk \oo{G}$?
\end{enumerate}

\section*{Acknowledgments}

The second author acknowledges support from the European Regional
Development Fund through the programme COMPETE
and by the Portuguese Government through FCT (Funda\c c\~ao para a Ci\^encia e a
Tecnologia) under the project PEst-C/MAT/UI0144/2011.


\clearpage
\addcontentsline{toc}{section}{References}
\begin{thebibliography}{99}
\bibitem{BZ} A.~Bjorner and G.~M.~Ziegler. Introduction to greedoids,
  in: {\em Matroid Applications} (ed. N. White), Cambridge
Univ. Press, pp. 284--357, 1992.
\bibitem{BT} R.~Brijder and L.~Traldi, The adjacency matroid of a
  graph, arXiv:1107.5493, preprint, 2011.
\bibitem{Cam} P.~J.~Cameron, Chamber systems and buildings, {\em The
    Encyclopaedia of Design Theory}, May 30, 2003. 
\bibitem{CFS} C.~R.~J.~Clapham, A.~Flockhart and J.~Sheehan, Graphs
  without four-cycles, {\em J. Graph Th.} 13.1 (1989), 29--47. 
\bibitem{Cox} H.~S.~M.~Coxeter, Self-dual configurations and regular
  graphs, {\em Bull. Amer. Math. Soc.} 56 (1950) 413--455.
\bibitem{Die}
R.~Diestel, {\em Graph Theory}, Springer-Verlag, 2000.
\bibitem{GHK} G.~Gierz, K.~H.~Hofmann, K.~Keimel, J.~D.~Lawson,
M.~Mislove and D.~S.~Scott, {\em Continuous lattices and domains},
volume 93 of {\em Encyclopedia of Mathematics and its Applications},
Cambridge University Press, Cambridge, 2003.
\bibitem{Gra} G.~Gr\"atzer, {\em Lattice theory: foundation}, Springer
  Basel AG, 2011. 
\bibitem{Gru} B.~Gr\"unbaum, {\em Configurations of points and lines},
  Graduate Studies in Mathematics, vol. 103, American Mathematical
  Society, 2009. 
\bibitem{Izh}
Z.~Izhakian, The tropical rank of a tropical matrix, preprint,
arXiv:math.AC/060420, 2006. 
\bibitem{IR1}
Z.~Izhakian and J.~Rhodes, New representations of monoids and
generalizations, preprint, arXiv:1103.0503, 2011.
\bibitem{IR2}
Z.~Izhakian and J.~Rhodes, Boolean representations of matroids and
lattices, preprint, arXiv:1108.1473, 2011.
\bibitem{IR3}
Z.~Izhakian and J.~Rhodes, C-independence and c-rank of posets and
lattices, preprint, arXiv:1110.3553, 2011.
\bibitem{Mar} G.~Markowski, Primes, irreducibles and extremal
  lattices, {\em Order} 9 (1992), 265--290.
\bibitem{MMT} R.~N.~McKenzie, G.~F.~McNulty and W.~F.~Taylor, {\em
Algebras, lattices, varieties, Vol. 1}, The Wadsworth $\&$ Brooks/Cole
Mathematics Series, Wadsworth $\&$ Brooks/Cole Advanced Books $\&$
Software, Monterey, CA, 1987.
\bibitem{Oxl} J.~G.~Oxley, {\em Matroid Theory}, Oxford Science
  Publications, 1992.
\bibitem{Oxl2} J.~G.~Oxley, What is a matroid, In: {\em LSU Mathematics
  Electronic Preprint Series}, pages 179--218, 2003. 
\bibitem{RSil} J.~Rhodes and P.~V.~Silva plus possibly other authors,
  Hereditary collections having boolean representations, in preparation.
\bibitem{RS} J.~Rhodes and B.~Steinberg, {\em The q-theory of Finite
    Semigroups}, Springer Monographs in Mathematics, 2009. 
\bibitem{Whi} H.~Whitney, On the abstract properties of linear
  dependence, {\em American Journal of Mathematics} (The Johns
Hopkins University Press), 57(3) (1935), 509--533 (Reprinted in Kung
(1986), pp. 55--79). 
\bibitem{W5}
Wikipedia,
http://en.wikipedia.org/wiki/Desargues$\underline{\hspace{.2cm}}$graph.
\bibitem{W4}
Wikipedia,
http://en.wikipedia.org/wiki/Desargues'$\underline{\hspace{.2cm}}$theorem. 
\bibitem{W7}
Wikipedia, http://en.wikipedia.org/wiki/Fano$\underline{\hspace{.2cm}}$plane.
\bibitem{W1}
Wikipedia, http://en.wikipedia.org/wiki/Heawood$\underline{\hspace{.2cm}}$graph.
\bibitem{W2}
Wikipedia, http://en.wikipedia.org/wiki/McGee$\underline{\hspace{.2cm}}$graph.
\bibitem{W8}
Wikipedia, http://en.wikipedia.org/wiki/Moore$\underline{\hspace{.2cm}}$graph.
\bibitem{W6}
Wikipedia, http://en.wikipedia.org/wiki/Table$
\underline{\hspace{.2cm}}$of$\underline{\hspace{.2cm}}$simple$
\underline{\hspace{.2cm}}$cubic$  
\underline{\hspace{.2cm}}$graphs.
\bibitem{W3}
Wikipedia, 
http://en.wikipedia.org/wiki/Tutte-Coxeter$\underline{\hspace{.2cm}}$graph. 
\end{thebibliography}
\end{document}